\documentclass{amsart}

\usepackage{graphicx,epsfig,amsfonts,amsmath}

\begin{document}


\def \build#1#2#3{\mathrel{\mathop{\kern 0pt#1}\limits_{#2}^{#3}}}
\def \B{\Big}
\def \b{\big}
\def \l{\left}
\def \r{\right}

\newdimen\AAdi%
\newbox\AAbo%
\font\AAFf=cmex10
\def\AArm{\fam0 \rm}%
\def\AAk#1#2{\setbox\AAbo=\hbox{#2}\AAdi=\wd\AAbo\kern#1\AAdi{}}%
\def\AAr#1#2#3{\setbox\AAbo=\hbox{#2}\AAdi=\ht\AAbo\raise#1\AAdi\hbox{#3}}%
\def\BBone{{\AArm 1\AAk{-.8}{I}I}}%
\def \build#1#2#3{\mathrel{\mathop{\kern 0pt#1}\limits_{#2}^{#3}}}

\newcommand{\rep}[1]{#1}

\newcommand{\bigO}[1]{\mathcal O\left( #1\right)}
\newcommand{\parth}[1]{\left(#1\right)}
\newcommand{\croc}[1]{\left[#1\right]}

\newcommand{\myepsfbox}[1]{
\epsfbox{\rep #1}
}
\newcommand{\mymarginpar}[1]{\marginpar{\raggedright{\tiny#1}}}
\newcommand{\PC}[1]{\mymarginpar{PC:\;#1}}

\newtheorem{theo}{Theorem}
\newtheorem{cor}{Corollary}
\newtheorem{proper}{Property}
\newtheorem{pro}{Proposition}
\newtheorem{lem}{Lemma}

\title[Random maps and ISE]{Random Planar Lattices and
Integrated SuperBrownian Excursion}

\author[Chassaing]{Philippe Chassaing}
\address{Philippe Chassaing,
Universit\'e Henri Poincar\'e,
B.P. 239,
54506 Vand{\oe}uvre-l\`es-Nancy}
\email{chassain@iecn.u-nancy.fr}

\author[Schaeffer]{Gilles Schaeffer}
\address{Gilles Schaeffer,
CNRS -- LORIA, B.P. 239,
54506 Vand{\oe}uvre-l\`es-Nancy}
\email{Gilles.Schaeffer@loria.fr, http://www.loria.fr/\string~schaeffe}

\begin{abstract}
In this paper, a surprising connection is described
between a specific brand of random lattices, namely planar
quadrangulations, and Aldous' Integrated SuperBrownian Excursion
(ISE).  As a consequence, the radius $r_n$ of a random quadrangulation
with $n$ faces is shown to converge, up to scaling, to the width
$r=R-L$ of the support of the one-dimensional ISE, or precisely:
\[
n^{-1/4}r_n\;\mathop{\longrightarrow}^{\textrm{\emph{law}}}\;(8/9)^{1/4}\,r.
\]
More generally the distribution of distances to a random vertex in a
random quadrangulation is
described in its scaled limit by the random measure ISE shifted to set
the minimum of its support in zero.

The first combinatorial ingredient is an encoding of quadrangulations
by \emph{trees embedded in the positive half-line}, reminiscent of
Cori and Vauquelin's well labelled trees.
The second step relates
these trees to embedded (discrete) trees in the sense of Aldous,
via the \emph{conjugation of tree principle}, an analogue for trees of
Vervaat's construction of the Brownian excursion from the bridge.

{From} probability theory, we need a new result of independent interest:
the weak convergence of the encoding of a random embedded
plane tree by two contour walks $(e^{(n)},\hat W^{(n)})$ to the
Brownian snake description $(e,\hat W)$ of ISE.

Our results suggest the existence of a \emph{Continuum Random Map}
describing in term of ISE the scaled limit of the dynamical
triangulations considered in two-dimensional pure quantum gravity.

\index{random planar maps}
\index{random lattices}
\index{ISE}
\index{CRT}
\index{discrete random trees}
\index{superprocesses}
\index{Brownian snake}
\end{abstract}

\maketitle

\section{Introduction}

{From} a distant perspective, this article uncovers a surprising, and
hopefully deep, relation between two famous models: \emph{random
planar maps}, as studied in combinatorics and quantum
physics, and \emph{Brownian snakes}, as studied in probability theory and
     statistical physics.  More precisely, our  results
connect some distance-related functionals of \emph{random
quadrangulations} with functionals of Aldous' \emph{Integrated
SuperBrownian Excursion} (ISE) in dimension one.

\subsection*{Quadrangulations} 
\begin{figure}
         \begin{center}
\includegraphics[width=11cm]{\rep 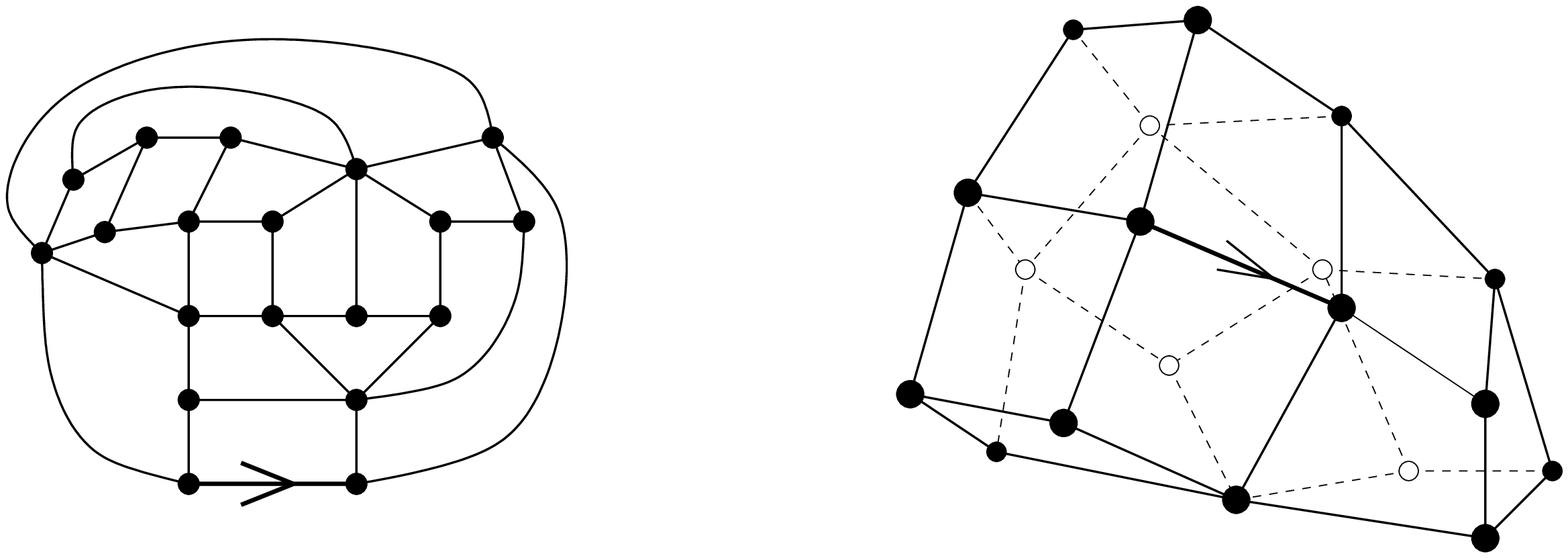}
         \end{center}
         \caption{Random quadrangulations, in planar or spherical
           representation.}
         \label{fig:aquad}
\end{figure}

On
the one hand, quadrangulations are finite plane graphs with
four-regular faces (see Figure~\ref{fig:aquad} and
Section~\ref{sec:def} for precise definitions). Random
quadrangulations, like random triangulations, random polyhedra, or the
$\phi^4$-models of physics, are instances of a general family of
random lattices that has received considerable attention both in
combinatorics (under the name \emph{random planar maps}, following
Tutte's terminology \cite{tutte-base}) and in physics (under the name
\emph{Euclidean two-dimensional discretised quantum geometry}, or
simply \emph{dynamical triangulations} or \emph{fluid lattices}
\cite{ambjorn, BIPZ, gross}).

Many probabilistic properties of random planar maps have been studied,
that are \emph{local properties} like vertex or face degrees
\cite{bender-canfield,gao-wormald}, or $0-1$ laws for properties expressible
in first order logic \cite{bender-compton}. Other well documented
families of properties are related to connectedness and constant size
separators \cite{BaFlScSo01}, also known as branchings into baby
universes \cite{phys-baby}. In this article we consider another
fundamental aspect of the geometry of random maps, namely \emph{global
properties of distances}. The \emph{profile}
$\smash{(H^{n}_k)_{k\geq0}}$ and \emph{radius} $r_n$ of a random
quadrangulation with $n$ faces are defined in analogy with the
classical profile and height of trees: $\smash{H_k^{n}}$ is the number
of vertices at distance $k$ from a basepoint, while $r_n$ is the maximal
distance reached. The profile was studied (with triangulations instead
of quadrangulations) by physicists Watabiki, Ambj{\o}rn et \emph{al.}
\cite{ambjorn-watabiki,watabiki} who gave a consistency argument proving
that the only possible scaling for the profile is $\smash{k\sim
n^{1/4}}$, a property which reads in their terminology \emph{the
internal Hausdorff dimension is 4}. Independently the conjecture that
$\mathbb{E}(r_n)\sim cn^{1/4}$ was proposed by Schaeffer
\cite{these-schaeffer}.

\subsection*{Integrated SuperBrownian Excursion}
On the other hand, ISE was introduced by Aldous as a model of random
distributions of masses \cite{AALD}.  He considers random embedded
discrete trees as obtained by the following two steps: first an
abstract tree $t$, say a Cayley tree with $n$ nodes, is taken from the
uniform distribution and each edge of $t$ is given length $1$; then
$t$ is embedded in the regular lattice on $\mathbb{Z}^d$, with the
root at the origin, and edges of the tree randomly mapped on edges of
the lattice. Assigning masses to leaves of the tree $t$ yield a random
distribution of mass on $\mathbb{Z}^d$. Upon scaling the lattice to
$n^{-1/4}\mathbb{Z}^d$, these random distributions of mass admit, for
$n$ going to infinity, a continuum limit $\mathcal{J}$ which is a
random probability measure on $\mathbb{R}^d$ called ISE.

Derbez and Slade proved that ISE describes in dimension larger than
eight the continuum limit of a model of lattice trees \cite{DS}, while
Hara and Slade obtained the same limit for the incipient infinite
cluster in percolation in dimension larger than six \cite{HS}. As
opposed to these works, we shall consider ISE in dimension one and our
embedded discrete trees should be thought of as folded on a
line.  The support of ISE is then a random interval $(L,R)$ of
$\mathbb{R}$ that contains the origin.

\subsection*{From quadrangulations to ISE}
The purpose of this paper is to draw a relation between, on the one
hand, random quadrangulations, and, on the other hand, Aldous' ISE:
upon proper scaling, the profile of a random quadrangulations is
described in the limit by ISE translated to have support $(0,R-L)$.
This relation implies in particular that the radius $r_n$ of random
quadrangulations, again upon scaling, weakly converges to the width
of the support of ISE in one dimension, that is the continuous random
variable $r=R-L$. We shall indeed prove
(Corollary~\ref{radcv}) that
\[
n^{-1/4}r_n\;\mathop{\longrightarrow}^{\textrm{law}}\;(8/9)^{1/4}\,r,
\]
as well as the convergence of moments. While this proves the
conjecture $\mathbb{E}(r_n)\sim cn^{1/4}$, the value of the
constant~$c$ remains unknown because, as mentioned by Aldous
\cite{AALD}, little is known on $R$ or $R-L$.
\medskip

The path from quadrangulations to ISE consists of three main steps,
the first two of combinatorial nature and the last with a more
probabilistic flavor. Our first step, Theorem~\ref{thm:welllab},
revisits a correspondence of Cori and Vauquelin \cite{cori-vauquelin}
between planar maps and some \emph{well labelled trees}, that can be
viewed as plane trees embedded in the positive half-line.  Thanks to
an alternative construction \cite[Ch.~7]{these-schaeffer}, we
show that under this correspondence the profile can be mapped to the
mass distribution on the half-line. In particular, the radius $r_n$
of a random quadrangulation is equal in law to the maximal label
$\mu_n$ of a random well labelled tree.

Safe for the positivity condition, well labelled trees would be
constructed exactly according to Aldous' prescription for embedded
discrete trees.  Well labelled trees are thus to Aldous' embedded trees
what the Brownian excursion is to the Brownian bridge, and we seek an
analogue of Vervaat's relation. At the discrete level a classical
elegant explanation of such relations is based on Dvoretsky and
Motzkin's cyclic shifts and cycle lemma. Our second combinatorial
step, Theorem~\ref{thm:coupling}, consists in the adaptation of these
ideas to embedded trees. More precisely, via the \emph{conjugation of
tree principle} of \cite[Chap. 2]{these-schaeffer}, we bound the
discrepancy between the mass distribution of our conditioned trees on
the positive half-line and a translated mass distribution of freely
embedded trees.  In particular we construct a coupling between well
labelled trees and freely embedded trees such that the largest label
$\mu_n$, and thus the radius $r_n$, is coupled to the width of the
support $(L_n,R_n)$ of random freely embedded trees:
\[
|r_n-(R_n-L_n)|\leq 3.
\]
Since our freely embedded trees are constructed according to Aldous'
prescription, one could expect to be able to conclude
directly. However two obstacles still need to be bypassed at this
point.

\subsection*{Contour walks and Brownian snakes.}
The first obstacle is that the construction of ISE as a continuum limit
of mass distributions supported by embedded discrete trees was only
outlined in Aldous' original paper. The original mathematical
definition is by embedding a continuum random tree (CRT), which
amounts to exchanging the embedding and the continuum limit.  But
Borgs et \emph{al.} proved that indeed ISE is the limit of mass
distributions supported by embedded Cayley trees \cite{BCHS} and their
proof could certainly be adapted to other simple classes of trees and
in particular to our embedded plane trees.

The second, more important, obstacle is that weak convergence
of probability measures is not adequate to our purpose, since we are
interested in particular in convergence of the width of the support,
\textit{which is not a continuous functional on the space of
measures}.  In order to circumvent this difficulty, we turn to the
description of ISE in terms of superprocesses: ISE can be constructed
from the Brownian snake with lifetime $e$, the standard Brownian
excursion \cite{AALD, LEG}.

{From} the discrete point of view, we consider the encoding of an
embedded plane tree by a pair of contour walks $(x_k,y_k)$, that
encode respectively the height of the node visited at time $k$ and its
position on the line.  Our last result,
Theorem~\ref{thm:converge}, is the weak convergence, upon proper
scaling, of this pair of walks to the Brownian snake with lifetime
$e$:
\[
\big(e^{(n)}(s),\hat W^{(n)}(s)\big)\mathop{\longrightarrow}^{\textrm{law}}
\big(e(s),\hat W_s\big).
\]
As $R=\sup_s \hat W_s$ and $L=\inf_s\hat W_s$ this convergence,
together with some deviation bounds obtained in the proof allows
us to conclude on the radius. (A similar weak convergence was
independently proved by Marckert and Mokkadem \cite{MM} but without
the deviation bounds we need here.)

More generally the joint convergence of the minimum and the mass
distribution of discrete embedded trees implies that, upon scaling, the
label distribution of well labelled trees converges to ISE translated
to have the minimum of its support at the origin. The same then holds
for the profile of random quadrangulations.

\subsection*{Dynamical triangulations and a Continuum Random Map}
Although we concentrate in this article on the radius and profile of
random quadrangulations, our derivation suggests a much tighter link
between random quadrangulations and ISE. We conjecture that a
Continuum Random Map (CRM) can be built from ISE that would describe
the continuum limit of scaled random quadrangulations, in a similar
way as the CRT describes the continuum limit of scaled random discrete
trees. From the point of view of physics, the resulting CRM would
describe in the limit the geometry of scaled dynamical triangulations
as studied in discretised two-dimensional Euclidean pure quantum
geometries \cite{ambjorn, BIPZ, gross}.  We plan to discuss this
connection further in future work.

\subsection*{Organization of the paper} Section~2 contains
the definition the combinatorial model of random lattice. Sections~3
and~4 are devoted to the first combinatorial steps. Section~5 
contains to the definition of probabilistic models, and the statement
of the convergence result. Finally in Section~6 we give the proof of this
convergence.

\begin{figure}
          \begin{center}
\includegraphics[width=11cm]{\rep 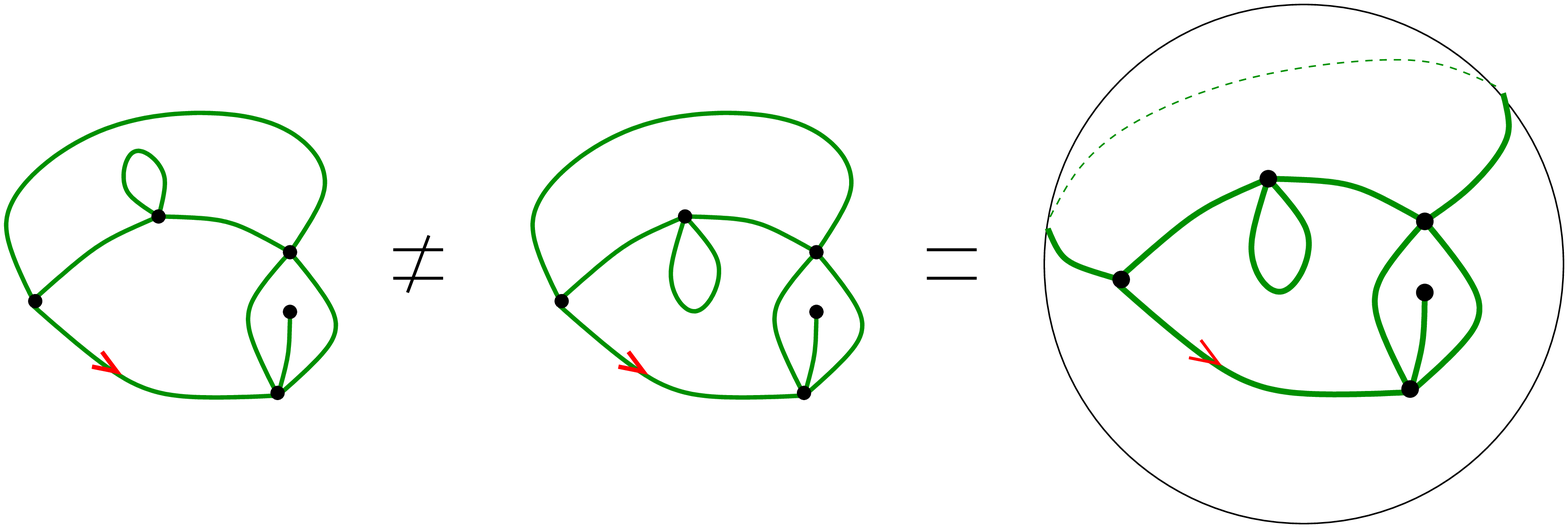}
          \end{center}
          \caption{Two distinct planar maps, and a spherical
representation of the
            second.}
          \label{fig:amap}
\end{figure}
\section{The combinatorial models of random lattice}\label{sec:def}

\subsection{Planar maps and quadrangulations}
A \textit{planar map} is a proper embedding (without edge crossings)
of a connected graph in the plane. Loops and multiple edges are
\emph{a priori} allowed. A planar map is \emph{rooted} if there is a
\emph{root}, \emph{i.e.}  a distinguished edge on the border of the
infinite face, which is oriented counterclockwise. The origin of the
root is called the \emph{root vertex}. Two rooted planar maps are
considered identical if there exists an homeomorphism \emph{of the
plane} that sends one map onto the other (roots included).

The difference between planar graphs and planar maps is that the
cyclic order of edges around vertices matters in maps, as illustrated
by Figure~\ref{fig:amap}.  Observe that planar maps can be equivalently
defined on the sphere.  In particular Euler's characteristic formula
applies and provides a relation between the numbers $n$ of edges, $f$
of faces and $v$ of vertices of any planar map: $f+v=n+2$.

\begin{figure}
         \begin{center}
\includegraphics[width=6cm]{\rep 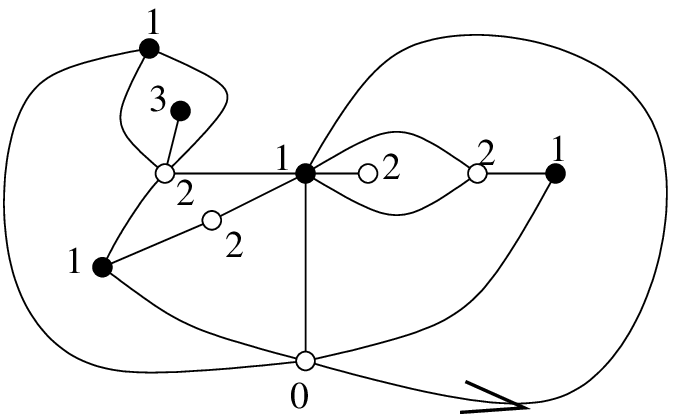}
\hspace{10mm}
\includegraphics[width=2cm]{\rep 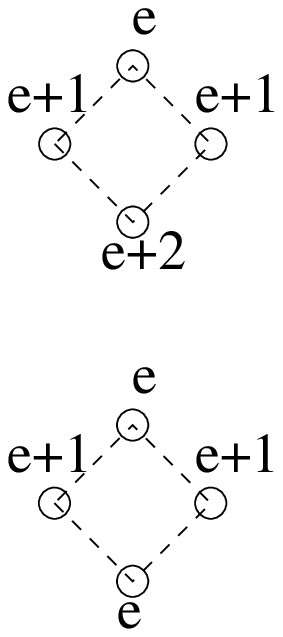}
         \end{center}
         \caption{Labelling by distance from the root vertex and the 
two possible
           configurations
           of labels (top: a simple face; bottom: a confluent face).}
         \label{fig:etiq}
\end{figure}

The \emph{degree} of a face or of a vertex of a map is its number of
incidence of edges. A planar map is a \emph{quadrangulation} if all
faces have degree four. All (planar) quadrangulations are
\emph{bipartite}: their vertices can be colored in black or white so
that the root is white and any edge joins two vertices with different
colors. In particular a quadrangulation contains no loop but may
contain multiple edges. See
Figures~\ref{fig:aquad} and~\ref{fig:etiq} for examples of quadrangulations.

Let $\mathcal{Q}_n$ denote the set of rooted quadrangulations with $n$
faces.  A quadrangulation with $n$ faces has $2n$ edges (because of
the degree constraint) and $n+2$ vertices (applying Euler's formula).
The number of rooted quadrangulations with $n$ faces was obtained by
W.T. Tutte \cite{tutte-base}:
\begin{equation}\label{equ:tutte}
|\mathcal{Q}_n|\;\;=\;\;\frac{2}{n+2}\,\frac{3^n}{n+1}{2n\choose n}.
\end{equation}
Various alternative proofs of this result have been obtained (see
\emph{e.g.} \cite{BIPZ, cori-vauquelin, arques-tres-bien-etiquete,
these-schaeffer}).  Our treatment will indirectly provide another
proof, related to \cite{cori-vauquelin, these-schaeffer}.

\subsection{Random planar lattices}
Let $L_n$ be a random variable with uniform distribution on
$\mathcal{Q}_n$. Formally, $L_n$ is the $\mathcal{Q}_n$-valued
random variable such that  for all
$Q\in\mathcal{Q}_n$
\[\Pr(L_n=Q)\;=\;\frac1{|\mathcal{Q}_n|}\;=\;
\frac1{\frac{2}{n+2}\frac{3^n}{n+1}{2n\choose
n}}.\]
The random variable $L_n$ is our \emph{random planar lattice}. To
explain this terminology, taken from physics, observe that locally the
usual planar square lattice is a planar map whose faces and vertices
all have degree four. Our random planar lattice corresponds to a
relaxation of the constraint on vertices.

Classical variants of this definition are obtained by replacing
quadrangulations with $n$ faces by triangulations with $2n$ triangles,
or by (vertex-)4-regular maps with $n$ vertices, or by all planar maps
with $n$ edges, \emph{etc.} All these random planar lattices have been
considered both in combinatorics (see \cite{BaFlScSo01} and references
therein) and in mathematical physics (see \cite{ambjorn} and
references therein; in the physics literature, definitions are
usually phrased using ``symmetry weights'' instead of rooted objects,
but this is strictly equivalent to the combinatorial
definition). Although details of local topology vary between families,
most probabilistic properties are believed to be ``universal'', that
is qualitatively analogue for all ``reasonable'' families. Observe also
that random maps in classical families have exponentially small
probability to be symmetric, so that all results hold as well as in
the model of uniform unrooted maps \cite{RiWo95}.

In this article we focus on quadrangulations because of their
combinatorial relation, detailed in Section~\ref{sec:welllab}, to well
labelled trees.

\subsection{The profile of a map}
The distance $d(x,y)$ between two vertices $x$ and $y$ of a map is the
minimal number of edges on a path from $x$ to $y$ (in other terms all
edges have abstract length $1$).

The \emph{profile} of a rooted map $M$ is the sequence
$(H_k)_{k\geq1}$, where $H_k\equiv H_k^{[M]}$ is the number of
vertices at distance $k$ of the root vertex $v_0$. We shall also
consider the cumulated profile $\widehat H_k^{[M]}\;=\;\sum_{\ell=1}^k
H_\ell^{[M]}$.  By construction the support of the profile of a rooted
map is an interval \emph{i.e.}  $\{k\mid H_k>0\}=[1,r]$ where $r$ is
the \emph{radius} of the map (sometimes also called
\emph{eccentricity}). The radius $r$ is closely related to the
\emph{diameter}, that is the largest distance between two vertices of
a map: in particular $r\leq d\leq 2r$.  The quadrangulation of
Figure~\ref{fig:etiq} has radius 3.

The \emph{profile of the random planar lattice $L_n$} is the random
variable $(H^{(n)}_k)_{k\geq1}$ that is defined by taking the profile
$(H^{[L_n]}_k)_{k\geq1}$ of an instance of $L_n$, while
     $(\widehat H^{(n)}_k)_{k\geq1}$  denotes the \emph{cumulated}
profile of $L_n$.
Similarly the radius of a random planar lattice
is a positive integer valued random variable $r_n$.

\section{Encoding the profile with well labelled trees}\label{sec:welllab}

\subsection{Well labelled trees and the encoding result}

\begin{figure}
\begin{center}
\includegraphics[width=6cm]{\rep 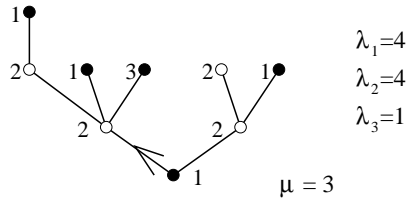}
\end{center}
\caption{A well labelled tree with its label distribution.
\label{fig:constrained}}
\end{figure}

A \emph{plane tree} is a rooted planar map without cycle (and thus
with only one face).  Equivalently plane trees can be recursively
defined as follows:
\begin{itemize}
\item the smallest tree is made of a single vertex,
\item any other tree is  a non-empty sequence
of subtrees attached to a root.
\end{itemize}
In other term,
each vertex has a possibly empty sequence of sons, and
each vertex but the root has a father.
The number of plane trees with $n$ edges is the well known
Catalan number
\[
C(2n)\;=\;\frac1{n+1}{2n\choose n}.
\]

A plane tree is \textit{well labelled} if all its vertices have
positive integral labels, the labels of two adjacent vertices differ
at most by one, and the label of the root vertex is one. Let
$\mathcal{W}_n$ denote the set of well labelled trees with $n$ edges.

The \emph{label distribution} of a well labelled tree $T$ is the
sequence $(\lambda_k)_{k\geq1}\equiv(\lambda_k^{[T]})_{k\geq1}$ where
$\lambda_k^{[T]}$ is the number of vertices with label $k$ in the tree
$T$.  The cumulated label distribution is defined by
$\widehat\lambda_k^{[T]}= \sum_{\ell=1}^{k}\lambda_\ell^{[T]}$. By
construction the support of the label distribution is an interval:
there exists an integer $\mu$ such that
$\{k\mid\lambda_k>0\}=[1,\mu]$. This integer $\mu$ is the maximal
label of the tree. These definitions are illustrated by
Figure~\ref{fig:constrained}.
\medskip

The following theorem will serve us to reduce the study of the profile
of quadrangulations to the study of the label distribution of well
labelled trees.
\begin{theo}[Schaeffer \cite{these-schaeffer}]\label{thm:welllab}
There exists a bijection $\mathcal{T}$ between rooted quadrangulations
with $n$ faces and well labelled trees with $n$ edges, such that the
profile $(H_k^{[Q]})_{k\geq1}$ of a quadrangulation $Q$ is mapped onto
the label distribution $(\lambda_k^{[T]})_{k\geq1}$ of the tree
$T=\mathcal{T}(Q)$.
\end{theo}
Theorem~\ref{thm:welllab} and Tutte's formula~(\ref{equ:tutte}) imply
that the number of well labelled trees with $n$ edges equals
\begin{equation}\label{equ:welllab}
|\mathcal{W}_n|\;=\;\frac{2}{n+2}\,\frac{3^n}{n+1}{2n\choose n}.
\end{equation}
This result was proved already by Cori and Vauquelin
\cite{cori-vauquelin}, who introduced well labelled trees to give an
encoding of all planar maps with $n$ edges.  Because of a classical
bijection between the latter maps and quadrangulations with $n$ faces,
their result is equivalent to the first part of Theorem~\ref{thm:welllab}.
Their bijection has been
extended to bipartite maps by Arqu\`es
\cite{arques-tres-bien-etiquete} and to higher genus maps by Marcus
and Vauquelin \cite{marcus-vauquelin}.  All these constructions were
recursive and based on encodings of maps with permutations (also known
as rotation systems).

However, our interest in well labelled trees lies in the relation
between the profile and the label distribution, which does not appear
in Cori and Vauquelin's bijection.  The bijection we use here is much
simpler and immediately leads to the second part of
Theorem~\ref{thm:welllab}. This approach was extended to non separable
maps by Jacquard \cite{benj-these} and to higher genus by Marcus and
Schaeffer \cite{marcus-schaeffer}.
\medskip

We postpone to Section~\ref{sec:conjugacy} the discussion of the
interesting form of Formula~(\ref{equ:welllab}) and its relation to
Catalan's numbers. Instead the rest of this part is concerned with the
proof of Theorem~\ref{thm:welllab}, which goes in three steps. First some
properties of distances in quadrangulations are indicated
(Section~\ref{sec:distances}). This allows in a second step to define
the encoding, as a mapping $\mathcal{T}$ from quadrangulations to well
labelled trees (Section~\ref{sec:encode}). A decoding procedure allows
then to prove that $\mathcal{T}$ is faithful
(Section~\ref{sec:decode}).

\subsection{Properties of distances in a quadrangulation}\label{sec:distances}

\begin{figure}
         \begin{center}
\includegraphics[width=2cm]{\rep 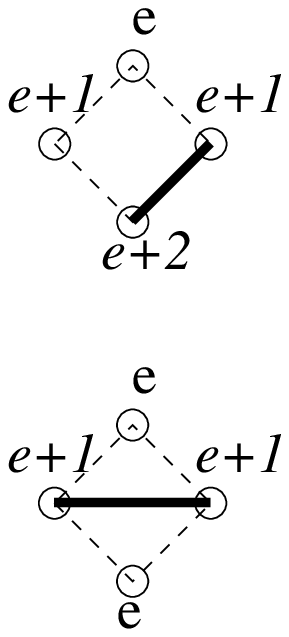}\hspace{3mm}
        \includegraphics[width=9cm]{\rep 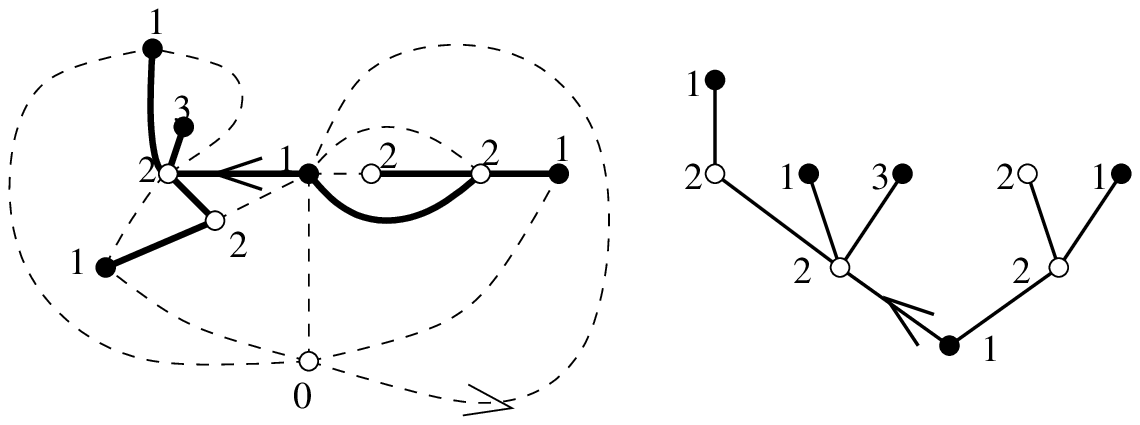}
         \end{center}
         \caption{The rules of selection of edges and an example.}
         \label{exemple}
\end{figure}

Let $Q$ be a rooted quadrangulation and denote $v_0$ its root vertex.
The labelling $\phi$ of the map ${Q}$ is defined by $\phi(x)=d(x,v_0)$
for each vertex $x$, where $d(x,y)$ denote the distance in $Q$
(cf. Figure~\ref{fig:etiq}). Observe that in the number of label $k$
in the labelling of the map $Q$ is precisely the number of vertices at
distance $k$ of $v_0$, that is $\smash{H^{[Q]}_k}=|\{x\mid
\phi(x)=k\}|$.
This labelling satisfies the following immediate properties:

\begin{proper} If $x$ and $y$ are joined by  an edge,
          $|\phi(x)-\phi(y)|=1$.  Indeed the quadrangulation being bipartite,
          a vertex $x$ is white if and only if  $\phi(x)$ is even, black if
          and only if  $\phi(x)$ is odd.
\end{proper}

\begin{proper} Around a face, four vertices appear: a black $x_1$,
          a white $y_1$, a black $x_2$ and a white
          $y_2$.  These vertices satisfy at least one of the two equalities
          $\phi(x_1)=\phi(x_2)$ or $\phi(y_1)=\phi(y_2)$ (cf.
          Figure~\ref{fig:etiq}).
\end{proper}

A face will be said \textit{simple} when only one equality is
satisfied and \textit{confluent} otherwise (see
Figure~\ref{fig:etiq}).  It should be noted that one may have
$x_1=x_2$ or $y_1=y_2$.

\subsection{Construction of the encoding $\mathcal{T}$}\label{sec:encode}
Let $Q$ be a rooted quadrangulation with its distance labelling. The
map ${Q}'$ is obtained by dividing all confluent faces ${Q}$ into two
triangular faces by an edge joining the two vertices with maximal
label.
Let us now define a subset $\mathcal{T}(Q)$ of edges of
${Q}'$ by two selection rules:
\begin{itemize}
\item In each confluent  face of ${Q}$, the edge that was
added to form ${Q}'$ is selected.
\item For each simple face $f$ of ${Q}$, an edge $e$ is selected: let
$v$ be the vertex with maximal label in $f$, then $e$ is the edge leaving
$v$ with $f$ on its left.
\end{itemize}
These two selection rules are illustrated by Figure~\ref{exemple}.
The first selected edge around the endpoint of the root of $Q$ is
taken to be the root of $\mathcal{T}(Q)$.

The proof of Theorem~\ref{thm:welllab} is now completed in two steps.
First, in the rest of this section, $\mathcal{T}(Q)$, which is \emph{a
priori} only defined as a subset of edges of ${Q}'$ together with
their incident vertices, is shown to be in fact a well labelled tree
with $n$ edges. Second, in the next section the inverse mapping is
described and used to prove that the mapping $\mathcal{T}$ is
faithful. \smallskip

\begin{pro}
The mapping ${\mathcal{T}}$ sends a quadrangulation $Q$ with $n$ faces
on a well labelled trees with $n$ edges.
\end{pro}
\begin{figure}
\begin{center}
\includegraphics[width=7cm]{\rep 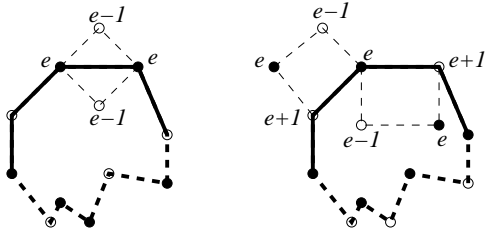}
\end{center}
\caption{Impossibility of cycles.}
\label{sans-cycle}
\end{figure}
\begin{proof}
If the vertex $x$ is not the root of $Q$, then one of its neighbors in
$Q$, say $y$, has a smaller label.  The edge $(x,y)$ can be incident
to: at least a confluent face; at least a simple face in which $x$ has
maximal label; or two simple faces in which $x$ has intermediate
label.  In all three cases, $x$ is incident to the selected edge of at
least one face. Thus all vertices of $Q$ but its root are also vertices of
$\mathcal{T}(Q)$: in particular $\mathcal{T}(Q)$ has $n+1$ vertices.
Next, the number of edges of $\mathcal{T}(Q)$ is $n$, because this is the
number of faces of $Q$ and two faces cannot select the same edge (as
immediately follows from inspection of the selection rules).  Now the
planarity of $Q$ and thus of ${Q}'$ grants that each connected
component of $\mathcal{T}(Q)$ is planar. Provided we can rule out
cycles, this imply that $\mathcal{T}(Q)$ is a forest of trees with $n$
edges and $n+1$ vertices, \emph{i.e.} a single tree. This tree is then
clearly well labelled.

Suppose now that there exists a cycle in $\mathcal{T}(Q)$ and
let $e\ge 0$ be
the value of the smallest label of a vertex of this cycle. Either all
these labels are equal to $e$, or there is in the cycle an edge
$(e,e+1)$ and an edge $(e+1,e)$. In both cases the rules of selection
of edges imply that each connected component of ${Q}'$ defined by the
cycle contains a vertex $x$ (resp. $y$) with label $e-1$, as shown by
Figure~\ref{sans-cycle}. According to Jordan's theorem,
either the shortest path from $x$ to the root
or the shortest path from $y$ to the root has to intersect the cycle,
leading to a contradiction
with the definition of labels by distances.  There are thus no cycles
and $\mathcal{T}(Q)$ is a tree.
\end{proof}

\subsection{The inverse $\mathcal{Q}$ of the mapping
$\mathcal{T}$}\label{sec:decode} Let $T$ be a well labelled tree with
$n$ edges. Recall that the tree $T$ can be viewed as a planar map that
has a unique face $F_0$. A \emph{corner} is a sector between two
consecutive edges around a vertex. A vertex of degree $k$ defines $k$
corners and the total number of corners of $T$ is $2n$. The label of a
corner is by definition the label of the corresponding vertex.

\begin{figure}
         \begin{center}
\includegraphics[width=11cm]{\rep 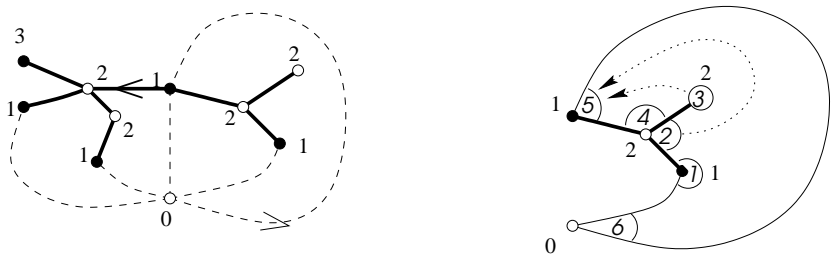}
         \end{center}
\caption{Step (1), and the cords $(i,s(i))$ in one of the faces.}
\label{fig:cords}
\end{figure}

The image $\mathcal{Q}(T)$ is defined in three steps.
\begin{enumerate}
\item A vertex $v_0$ with label $0$ is placed in the face $F_0$ and
one edge is added between this vertex and each of the $\ell$ corners
with label $1$. The new root is taken to be the edge arriving from
$v_0$ at the corner before the root of $T$.
\end{enumerate}
After Step (1) a uniquely defined rooted map $T_0$ with $\ell$ faces
has been obtained (see Figure~\ref{fig:cords}, with $\ell=5$). The
next steps take place independently in each of those faces and are
thus described for a generic\footnote{The infinite face is only
apparently different from the others: imagine the map on the sphere.}
face $F$ of $T_0$.  Let $k$ be the degree of $F$ (by construction
$k\geq3$).  Among the corners of $F$ only one belongs to $v_0$ and has
label $0$. Let the corners be numbered from $1$ to $k$ in clockwise
order along the border, starting right after $v_0$. Let moreover $e_i$
be the label of corner $i$ (so that $e_1=e_{k-1}=1$ and $e_k=0$). In
Figure~\ref{fig:cords} the corners are explicitly represented with
their numbering for one of the faces.

\begin{enumerate}
\setcounter{enumi}{1}
\item
The function successor $s$ is defined for all corners $1,\ldots,k-1$
by
\[
s(i)= \inf\{j>i\mid e_j=e_i-1\}.
\]
For each corner $i\ge 2$ such that $s(i)\neq i+1$, a cord $(i,s(i))$ is
added inside the face, in such a way that the various cords do not
intersect (Property~\ref{pro:licite}).
\end{enumerate}
Once this construction has be carried on in each face, a planar map
${T}'$ is obtained.
\begin{enumerate}
\setcounter{enumi}{2}
\item All edges of with labels of the form $(e,e)$ of ${T}'$ are
deleted. The resulting map is a quadrangulation $\mathcal{Q}(T)$
with $n$ faces (Property~\ref{pro:isquad}).
\end{enumerate}
The following proposition ends the proof of Theorem~\ref{thm:welllab}.
\begin{pro}\label{pro:inverse}
The mapping $\mathcal Q$ is the inverse bijection of the mapping
$\mathcal{T}$.
\end{pro}
Let us first prove the two properties that validate the preceding
construction.
\begin{proper}\label{pro:licite}
The cords $(i,s(i))$ do not intersect.
\end{proper}
\begin{proof}
Suppose that two cords $(i,s(i))$ and $(j,s(j))$ cross each
other.  Upon maybe exchanging $i$ and $j$ one has $i<j<s(i)<s(j)$.
The first two inequalities imply, together with the definition of $s$,
that $e_j> e_{s(i)}$, while the two last inequalities imply
$e_{s(i)}\ge e_j$. This is a contradiction.
\end{proof}

\begin{proper}\label{pro:isquad}
The faces of ${T}'$ are of one of the two types of
Figure~\ref{fig:deuxfaces}: either triangular with labels $e,e+1,e+1$,
or quadrangular with labels $e,e+1,e+2,e+1$. The faces of
$\mathcal{Q}(T)$  are all quadrangular.
\end{proper}
\begin{proof} Let $f$ be a face of ${T}'$. The face $f$ is included in
a face $F$ of $T_0$ so that its corners inherit the numbering and
labelling of those of $F$.  Let $j$ be the corner with largest number
in $f$ and $i_1<i_2<j$ its two neighbors in $f$ (cf.
Figure~\ref{fig:deuxfaces}).
Let us compute the label of the corners $i_1$ and $i_2$:
the edge  $(i_1,j)$ is a cord by construction so that $j=s(i_1)$ and
$e_{i_1}=e_j+1$; moreover, as $i_1<i_2<j$, this imply $e_{i_2}\ge
e_{i_1}$ (or $j$ would not be $s(i_1)$) and finally $e_{i_2}=e_{i_1}=e_j+1$.

By construction and planarity, no cord can arrive at $i_1$ between
the unique leaving cord $(i_1,j)$ and the edge $(i_1,i_1+1)$ of $F$.
The latter edge thus borders the face $f$. There are then two cases,
as illustrated by Figure~\ref{fig:deuxfaces}:
\begin{itemize}
\item if $e_{i_1+1}=e_{i_1}$, then $i_2=i_1+1$, and the face is
triangular (left hand figure),
\item otherwise $e_{i_1+1}=e_{i_1}+1$ (recall $e_{i_1+1}\geq e_{i_1}$
since $i_1<i_1+1<j=s(i_1)$) and the cord leaving $i_1+1$ goes to $i_2$
(otherwise $s(i_1+1)<i_2$ with $e_{s(i_1+1)}=e_{i_1+1}-1=e_{i_2}$ and
a cord $(s(i_1+1),j)$ would exclude $i_2$ from the face $f$): the
face is quadrangular (right hand side in Figure~\ref{fig:deuxfaces}).
\end{itemize}
Observe finally that the deletion of edges with labels of the form
$(e,e)$ will join triangular faces two by two so that $\mathcal{Q}(T)$
has only quadrangular faces.
\begin{figure}
\begin{center}
\includegraphics[width=12cm]{\rep 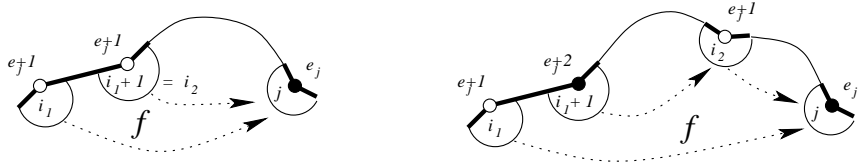}
\end{center}
\caption{Two possible sizes for $f$: triangular or
quadrangular.\label{fig:deuxfaces}}
\end{figure}
\end{proof}
\begin{proof}[Proof of Proposition~\ref{pro:inverse}] Given a well
labelled tree $T$, faces of its image $\mathcal{Q}(T)$ are as
described by Figure~\ref{fig:deuxfaces}. The selection rules for
$\mathcal{T}$ then shows that each face correctly selects an edge of
$T$, so that ${\mathcal T}({\mathcal Q}(T))=T$.  Thus ${\mathcal T}$
and ${\mathcal Q}$ are inverse bijections between well labelled trees
and a subset $\tilde{\mathcal Q}_n$ of the set of
quadrangulations. Equality of cardinalities, as granted for instance
by the alternative bijection of \cite{cori-vauquelin}, proves that
$\tilde {\mathcal Q}_n$ is the full set of quadrangulations with $n$
faces and concludes the proof.  However, we provide below a direct
proof of the equality ${\mathcal Q}({\mathcal{T}}(Q))=Q$ for any
quadrangulation $Q$, for this provides better understanding of the
bijection.

Let $Q$ be a quadrangulation, $Q'$ and $T=\mathcal{T}(Q)$ the map and
tree as in Section~\ref{sec:encode}, and $T_0$, $T'$ as in
the construction of $\mathcal{Q}(T)$.  Consider first the
selection rules applied around the root to construct $\mathcal
T(Q)$. Each edge (with labels) 1-2 of $\mathcal T(Q)$ forms a directly
oriented corner with an edge 1-0 in its face of creation, while each
edge 1-1 forms two such corners (one on each side). Hence, in
accordance with Step~(1) of the reciprocal construction, an edge 1-0
arrives at each corner with label $1$. Thus the submap $T_0$ of $T'$
is also a submap of $Q'$.
Moreover $T_0$ covers all vertices of $T'$
(resp. $Q'$), so that edges of $T'$ (resp. $Q'$) not in $T_0$ are
cords of faces of $T_0$. Accordingly, $T'=Q'$ if, inside each face of
$T_0$, both $T'$ and $Q'$ have the same cords.

The maps ${{Q}'}$ and ${T}'$ have the same vertices, and, due to
Property~\ref{pro:isquad}, the same number of faces of degree~4 (that
is, the number of edges $i-(i+1)$ in $T$), and the same
number of faces of degree 3 (that is, the number of edges $i-i$ in
$T$). By Euler's formula, they thus have also the same number of
edges and finally, they are equal as planar maps if, inside all faces
of $T_0$, each cord of ${T}'$ is a cord of ${{Q}'}$.

Let us now work inside a face $F$ of $T_0$ (see
Figure~\ref{fig:faceinduction}). By construction of $T_0$ the face $F$
has only one corner with label $0$ (incident to the root) and two
corners with label $1$ (since such a corner is incident to an edge
$(0,1)$ after Step~(1)). If $F$ has degree 3 (resp. 4), the corners
around the face, in clockwise order, are labelled 0-1-1
(resp. 0-1-2-1), and there is no cord, neither in $Q'$ nor in
$T'$. Let us thus assume that $F$ has degree $k$ larger than 4, and
number the $k$ corners of $F$ in clockwise order starting after the
root corner. Let $e_i$ be the label of the $i$th corner, so that
$e_1=1$, $e_2=2$, $e_{k-2}=2$, $e_{k-1}=1$, $e_k=0$, and $e_i\geq2$
otherwise. This numbering corresponds to the one used to define
cords of $T'$: for each corner $i$ with $e_i=2$ but
the last one, $s(i)=k-1$ and the cord $(i,k-1)$ appears in $T'$.

In order to check that these cords appear also in $Q'$ we consider
corners with label 2 in increasing order: they are numbered
$2=i_1<i_2<\ldots<i_p=k-2$. In $Q'$ let $f_1$ be the face that
contains the corners numbered $k$, $1$ and $k-1$ (with
labels $0$, $1$ and $1$). The face $f_1$ contains a fourth corner, with
label $2$: it can be $i_p$ (if the cord $(1,i_p)$ is in
$Q'$) or $i_1$ (if the cord $(i_1,k-1)$ is in $Q'$). In the
first case the face $f_1$ in $Q'$ has corners $(1,i_p,k-1,k)$ and there
is a contradiction with the fact that the edge $(i_p,k-1)$ has not
been selected in $F$ by the construction $\mathcal{T}$.  Hence the
second case hold: the cord $(i_1,k-1)$ is in $Q'$ and the edge
$(1,i_1)$ was selected.

\begin{figure}
\begin{center}
\includegraphics[width=4.2cm]{\rep 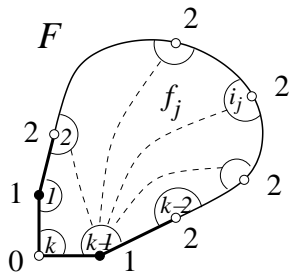}
\end{center}
\caption{Faces of $T_1$ inside a face  $F$ of $T_0$.\label{fig:faceinduction}}
\end{figure}

Now assume that the edges $(i_{j'},k-1)$ belong to $Q'$ for $j'<j<k-2$
and check that $(i_j,k-1)$ belongs to $Q'$. Consider in $Q'$ the face
$f_j$, included in the face $F$ and bordering the edge $(i_j-1,i_j)$
of the cycle $F$. If $e_{i_j-1}=2$ then $i_j-1=i_{j-1}$ and the face
$f_j$ is triangular (since the selection rule of confluent faces was
applied by $\mathcal{T}$) and contains the cord $(i_j,k-1)$.
Otherwise $e_{i_j-1}=3$ and the face $f_j$ is quadrangular and of
simple type (since the selection rule of simple faces was applied by
$\mathcal{T}$). Therefore there is an edge from $i_j$ to a corner with
label $1$, which can only be the cord $(i_j,k-1)$, for a cord
$(i_j,1)$ would cross $(i_1,k-1)$.

All cords with labels 1-2 are thus identical in $T'$ and $Q'$. Let
$T_1$ be the union of $T_0$ and these cords. In view of the previous
discussion, the faces of $T_1$ are exactly the previous subdivisions
into faces $f_j$ of all faces $F$ of $T_0$. Moreover each face $f$ of
$T_1$ has only one corner with label~1 and two with label~2, all other
labels being at least~3. Shifting down all labels by one inside face
$f$, the situation is exactly equivalent to that of the face $F$ above
(observe that the rules for the construction $\mathcal{T}$ and
$\mathcal{Q}$ remain unaffected by the shift since only vertices with
label greater or equal to two are considered). The identity between
cords of $T'$ and $Q'$ of successively greater orders can thus be
checked inductively.

Finally ${{Q}'}$ contains all the edges of ${T}'$, that is
     ${{Q}'}={T}'$, and since the deletion of edges with labels of the
     form $(e,e)$ in ${{Q}'}$ (resp. ${T}'$) produces $Q$
     (resp. ${\mathcal Q}({\mathcal{T}}(Q))$), we obtain
     ${\mathcal Q}({\mathcal{T}}(Q))=Q$.  \end{proof}

\section{Well labelled and embedded trees}\label{sec:conjugacy}

\subsection{Unconstrained well labelled trees as embedded trees}
Formula~(\ref{equ:welllab}) for the number of well labelled trees with
$n$ edges,
\[
|\mathcal{W}_n|\;=\;\frac{2}{n+2}\,\frac{3^n}{n+1}{2n\choose n}
\;=\;\frac{2}{n+2}\cdot 3^n\cdot C(2n),
\]
is remarkably simple and yet not immediately clear from definition.
Indeed, even though $C(2n)$ is known to be the number of plane trees,
the positivity of labels makes it difficult to count labellings that
make a plane tree well labelled.

It is thus natural to work first without this positivity condition:
define a plane tree to be an \textit{unconstrained well labelled tree}
if its vertices have integral labels, the labels of two adjacent
vertices differ at most by one, and the label of the root vertex is
one.  Let $\mathcal{U}_n$ denote the set of unconstrained well
labelled trees with $n$ edges.

The labelling of a labelled tree can be recovered uniquely from the
label of its root and the variations of labels along all edges. We
shall denote $\kappa(\epsilon)\in\{-1,0,1\}$ the variation of labels
along the edge $\epsilon$ when it is traversed away from the root.
Since there is no positivity condition on the labels of unconstrained
well labelled trees, all $\kappa(\epsilon)$ can be set independently
and the number of labellings of a plane tree that yield an
unconstrained well labelled tree is just $3^n$. That is,
\[
|\mathcal{U}_n|\;=\;\frac{3^n}{n+1}\binom{2n}{n}\;=\;3^n\cdot C(2n).
\]

\begin{figure}
\begin{center}
\includegraphics[width=12cm]{\rep 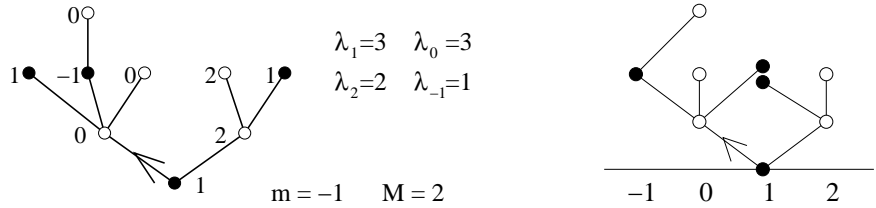}
\end{center}
\caption{An unconstrained well labelled tree with its label
distribution and a representation of the embedding on the line
(the plane order structure of the tree is lost in the latter representation).
\label{fig:unconstrained}}
\end{figure}

The definition of label distribution extends to unconstrained well
labelled trees. For $U\in\mathcal{U}_n$ let
$(\lambda_k)_{m<k<M}\equiv(\lambda_k^{[U]})_{k\in\mathbb{Z}}$ be the
number of vertices with label $k$ in the tree $U$.  The label
distribution of $U$ is supported by an interval $[m,M]$ with $m\leq
1\leq M$. The cumulated label distribution is defined with respect to
the minimum label $m$ by
$\widehat\lambda_k^{[U]}=\sum_{\ell=1}^{k}\lambda_{m+\ell-1}^{[U]}$.
These definitions are illustrated by Figure~\ref{fig:unconstrained}.

Observe moreover that similar unconstrained labellings have been
considered by D. Aldous \cite{AALD} with the following interpretation
(we restrict to our special one-dimensional case).
The tree is folded on the lattice $\mathbb{Z}$ with the root set at
position $1$ and each edge mapped on an elementary vector
(here $+1$, $0$, or $-1$). The label of a node then describe its
position on the line and, upon counting the number of nodes at
position $j$, a mass distribution is obtained. More precisely,
with our notations, Aldous' discrete mass distribution associated to a
tree $U\in \mathcal{U}_n$ is just the empirical measure of labels
\[
\mathcal{J}^{[U]}\;=\;\frac{1}n\sum_{k\in\mathbb{Z}}\lambda^{[U]}_k\delta_k,
\]
where $\delta_k$ denote the dirac mass at $k$.

In view of this interpretation and for concision's sake, let us rename
\emph{unconstrained well labelled trees} and call them instead
\emph{embedded trees}.

\subsection{Random trees and random quadrangulations}
Let $W_n$ and $U_n$ be random variables with uniform distribution on
$\mathcal{W}_n$ and $\mathcal{U}_n$. More precisely,
\[
\Pr(W_n=W)\;=\;\frac1{\frac{2}{n+2}\frac{3^n}{n+1}{2n\choose n}},
\quad\textrm{ and }\quad
\Pr(U_n=U)\;=\;\frac1{\frac{3^n}{n+1}{2n\choose n}},
\]
for all $W\in\mathcal{W}_n$ and $U\in\mathcal{U}_n$.

The label distribution of the corresponding random trees are two
random variables that we shall denote
${(\lambda^{(n)}_k)_{k\geq1}\equiv(\lambda^{[W_n]}_k)_{k\geq1}}$ for
random well labelled trees, and
${(\Lambda^{(n)}_k)_{k\in\mathbb{Z}}\equiv
(\lambda^{[U_n]}_k)_{k\in\mathbb{Z}}}$ for random embedded trees. For
random well labelled trees we also use the notation $\mu_n$ for the
maximal label, and for random embedded trees the notations $m_n$ and
$M_n$ for the minimal and maximal label respectively. Finally
cumulated profiles
$\widehat\lambda_k^{(n)}=\sum_{\ell=1}^k\lambda_\ell^{[W_n]}$ and
$\widehat\Lambda_k^{(n)}=\sum_{\ell=1}^k\lambda_{m_n+\ell-1}^{[U_n]}$ are
defined accordingly (the minimum $m_n$ in $\widehat\Lambda_k^{(n)}$ is
understood for the same realisation $U_n$).
\medskip

At this point we are given three random variables: random
quadrangulations $L_n$, random well labelled trees $W_n$ and random
embedded trees $U_n$. On the one hand, according to
Theorem~\ref{thm:welllab}, random quadrangulations ``are'' random well
labelled trees, as illustrated by the next corollary.
\begin{cor}
\label{cor:id}
The label distribution of random well labelled trees has the same
distribution as the profile of quadrangulations:
\[
(\lambda^{(n)}_k)_{k\geq1}\;\;\mathop{=\rule{0mm}{2mm}}^{\textrm{law}}\;\;
(H^{(n)}_k)_{k\geq1}.
\]
         In particular $r_n=\mu_n$.
\end{cor}

On the other hand, random embedded trees seem to be a simple variant
of well labelled trees that has the great advantage to be defined in
accordance with Aldous' prescription for discrete embedded trees.
This leads us to study more precisely the relation between $W_n$ and $U_n$.
By definition, $\mathcal{W}_n\subset\mathcal{U}_n$, and according to
Tutte's formula~(\ref{equ:welllab}),
\begin{equation}\label{equ:UtoW}
|\mathcal{W}_n|\;=\;\frac{2}{n+2}\cdot|\mathcal{U}_n|.
\end{equation}
For  combinatorists, this relation could be reminiscent of
the relation between the number of Dyck walks and the number of
bilatere Dyck walks (see \cite[Ch. 5]{Stanley}).

Equivalently, from a more probabilistic point of view, the relation reads
\[
\Pr(U_n\in\mathcal{W}_n)\;=\;\frac{2}{n+2},
\]
and random well labelled trees are random embedded trees conditioned
to positivity. This is exactly similar to Kemperman's formula for the
probability that a simple symmetric walk on $\mathbb{Z}$ starting from
$k>0$ and ending at $0$ after $n$ steps remains positive until the last step
(see \cite{pitman}).

\subsection{Cyclic shifts and the cycle lemma}\label{sec:cycle}
The idea to consider cyclic shifts originates in Dvoretsky and
Motzkin's work and was used by Raney to prove Lagrange inversion
formula and by Tak\'acs to prove and extend Kemperman's formula for
random walks (we refer to \cite[Ch. 5]{Stanley} and \cite{pitman} for
these historical references and many more). We shall prove a
consequence of this idea to the study of ``height distribution'' of
simple walks, that will be fundamental in the next section.

Let $n$ and $k$ be nonnegative integers and let $\mathcal{B}_{n,k}$ denote
the set of walks of length $2n+k$ with $n$ increments $+1$ and $n+k$
increments $-1$, that end with a negative increment $-1$. A walk
$w\in\mathcal{B}_{n,k}$ is described either by its sequence of
increments $w=(w_1,\ldots,w_{2n+k})$, $w_i\in\{+1,-1\}$, or by the
partial sums $w(p)=\sum_{i=1}^pw_i$, $p=0,\ldots,2n+k$.  By construction,
$w(0)=0$, $w_{2n+k}=-1$, $w(2n+k)=-k$, and $w(2n+k-1)=-k+1$.
Finally for $k\geq1$, consider the subset $\mathcal{D}_{n,k}$ of
$\mathcal{B}_{n,k}$ of ``positive'' walks defined by the condition:
$w(p)>-k$ for all $0\leq p<2n+k$.

Two walks $w$ and $w'$ of $\mathcal{B}_{n,k}$ belong the same
\emph{conjugacy class} if they differ by a cyclic shift, that is, if
there exists $s$ such that
\[
w=(w_1,\ldots,w_{2n+k})\quad\textrm{ and }\quad
w'=(w_s,\ldots,w_{s+2n+k}),
\]
where indices are considered modulo $2n+k$.

Define a (left-to-right) record to be a step $p\geq1$ at which a
minimum is reached for the first time: $w(q)>w(p)$ for all $q<p$.
Since $w(2n+k)=-k$ there are at least $k$ records. Let us denote
$p_1<\cdots<p_k$ the $k$ lowest records, so that in particular $w(p_k)$
is the minimum value reached by the walk and $w(p_i)=k-i+w(p_k)$.  The
steps $p_i$, $i=1,\ldots,k$ are called the \emph{low records} of $w$.

\begin{figure}
\begin{center}
\includegraphics[width=7cm]{\rep 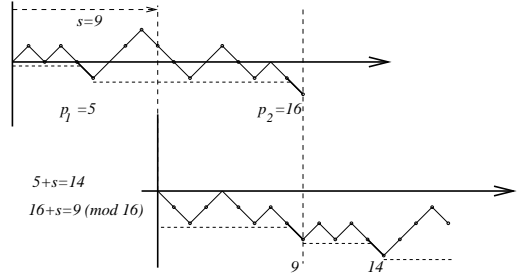}
\end{center}
\caption{Two conjugated walks from $\mathcal{D}_{n,2}$ and
$\mathcal{B}_{n,2}$ with low records and floor level of Dyck height indicated.
\label{fig:cyclic}}
\end{figure}

The following immediate properties are illustrated by Figure~\ref{fig:cyclic}.
\begin{proper}\label{pro:Dyck}\hspace{-1mm}
A walk of $\mathcal{B}_{n,k}\!$ belongs to $\mathcal{D}_{n,k}\!$ if and
only if its lowest record $p_k\!=\!2n\!+\!k$.
\end{proper}
\begin{proper}\label{pro:records}
Cyclic shifts transport low records:
Let
$w=(w_1,\ldots,w_{2n+k})$ and $w'=(w_s,\ldots,w_{s+2n+k})$ and
assume $\{p_1,\dots,p_k\}$ are the low records of
$w$. Then the low records of $w'$ are
$\{p_1+s,\ldots,p_k+s\}$ (modulo $2n+k$).
\end{proper}

The classical cycle lemma follows from
Properties~\ref{pro:Dyck} and~\ref{pro:records}.
\begin{lem}[Cycle lemma]\label{lem:cycle}
Let ${C}$ be a conjugacy class of $\mathcal{B}_{n,k}$. Then
\[
(n+k)\cdot|{C}\cap\mathcal{D}_{n,k}|\;=\;k\cdot|{C}|.
\]
\end{lem}
\begin{proof}
Let us apply a double counting argument:
\begin{itemize}
\item The left hand side counts walks in ${C}$ with a (low)
record at the last position ($w\in\mathcal{D}_{n,k}$) and a down step
marked ($n+k$ choices).
\item
The right hand side counts walks in ${C}$ with a down step
at the last position ($w\in\mathcal{B}_{n,k}$) and a low record marked
($k$ choices).
\end{itemize}
Now a bijection is obtained between these two sets upon sending the
marked step to the last position by a cyclic shift and marking the
former last step.
\end{proof}

Given a walk $w$ with lowest record $p_k$, the \emph{height-to-min} of a step
$p$ is $\tilde w(p)=w(p)-w(p_k)$, which is nonnegative by definition.
In order to study $\tilde w$, it will be convenient to consider the
height of the walk relatively to the $k$ low records: given a walk $w$
with low records $p_1<\cdots<p_k$, let us define the \emph{Dyck height} at
step $p$ by
\[
\bar w(p)=
\left\{
\begin{array}{ll}
w(p)-w(p_1)+1, & \textrm{if }0\leq p<p_1,\\
w(p)-w(p_i), & \textrm{if }p_i\leq p<p_{i+1}, \textrm{ with }1\leq i\leq k-1\\
w(p)-w(p_k), & \textrm{if }p_k\leq p\leq2n+k.
\end{array}
\right.
\]
The Dyck height can be understood as the height inside each of
the $k$ Dyck factors separated by low records.  Let $\hat\ell_i(w)$
(resp. $\hat h_i(w)$) denote the number of down steps of $w$ ending at
Dyck height (resp. height-to-min) at most $i$,
\[
\smash{\hat\ell_i}(w)\;=\;|\{p\mid \bar w(p)\leq i, \; w_p=-1\}|,\;\;
\textrm{ and }\;\;
\smash{\hat h_i}(w)\;=\;|\{p\mid \tilde w(p)\leq i, \; w_p=-1\}|.
\]
Then by construction, for all $w$ and $p$, $\bar w(p)\leq\tilde
w(p)\leq\bar w(p)+k$ and, for all $i$, $\hat
h_i(w)\leq\hat\ell_i(w)\leq\hat h_{i+k}(w)$.  The following lemma
immediately follows from Property~\ref{pro:records}.
\begin{lem}\label{lem:height}
The Dyck height commutes with cyclic shift. In particular
the Dyck height distribution $\hat\ell_i$ is invariant under cyclic shift:
\[
(\hat\ell_i(w))_{i\geq0}\;=\;(\hat\ell_i(w'))_{i\geq0},\qquad
\textrm{for all $w$ and $w'$ in the same conjugacy class.}
\]
The height-to-min distribution thus satisfies the following weaker
invariance:
\[
\hat h_i(w)\;\leq\; \hat h_{i+k}(w'), \qquad \textrm{for all $i\geq 0$
and $w$, $w'$ in the same conjugacy class.}
\]
\end{lem}

    From  the probabilistic point of view this result can be understood as
a simplified discrete version of Vervaat's relation between the
Brownian excursion and the Brownian bridge and their local times
relatively to the minimum.

\subsection{How to lift the positivity condition for labelled trees}
In view of Relation~(\ref{equ:UtoW}) one can expect to apply ideas of
the previous section to related well labelled trees to embedded trees.
As a matter of fact we shall prove the following theorem.
\begin{theo}\label{thm:combi}
There exists a partition of
$\mathcal{U}_n=\bigcup_{C\in\mathcal{C}_n}C$ into disjoint
\emph{conjugacy classes} each of size at most $n+2$ and such that
in each class $C\in\mathcal{C}_n$
\begin{itemize}
\item well labelled trees are fairly represented:
\[
2\cdot|{C}|\;=\;(n+2)\cdot|{C}\cap\mathcal{W}_n|,
\]
\item and for any $W\in\mathcal{W}_n\cap C$, $U\in C$ and $k\geq1$,
\[
\widehat\Lambda_{k-2}(U)\;\leq\;
\widehat\lambda_k(W)\;\leq\;
\widehat\Lambda_{k+2}(U).
\]
\end{itemize}
\end{theo}
\begin{cor}[Cori-Vauquelin, 1981]\label{cor:count}
The number of well labelled trees with $n$ edges, (which is also the
number of quadrangulations with $n$ faces), is
\[
|\mathcal{W}_n|\;=\;\frac{2}{n+2}\cdot|\mathcal{U}_n|
\;=\;\frac{2}{n+2}\cdot\frac{3^n}{n+1}\binom{2n}{n}.
\]
\end{cor}
The proof of Theorem~\ref{thm:combi} is presented in the next
section. It relies on an encoding of plane trees in terms of another
family of trees, called \emph{blossom trees}, and on the
\emph{conjugation of trees} principle which is an analogue of the cycle
lemma for blossom trees. This principle was introduced
in~\cite{these-schaeffer} in order to give a direct combinatorial
proof of Corollary~\ref{cor:count} based on the cycle lemma. However
that proof did not rely on well labelled trees and does not provide
the link to the profile.

Theorem~\ref{thm:combi} admits the following probabilistic restatement.
\begin{theo}\label{thm:coupling}
There is a coupling $(W_n,U_n)$ (\emph{i.e.} a distribution on
$\mathcal{W}_n\times\mathcal{U}_n$ such that the marginals are
$W_n$ and $U_n$ as previously defined)
such that the induced joint distribution $(\lambda^{(n)},\Lambda^{(n)})$
satisfies for all $k$
\[
\widehat\Lambda_{k-2}^{(n)}\;\leq\;
\widehat\lambda_k^{(n)}\;\leq\;
\widehat\Lambda_{k+2}^{(n)},
\]
and in particular
\[
|\mu_n-(M_n-m_n)|\leq3.
\]
\end{theo}
\begin{proof}[Proof of Theorem~\ref{thm:coupling}]
The distribution on $\mathcal{W}_n\times\mathcal{U}_n$ is immediately
obtained from the partition
$\mathcal{U}_n=\bigcup_{C\in\mathcal{C}_n}C$ as follows: for any $(W,U)$ in
$\mathcal{W}_n\times\mathcal{U}_n$, let
\[
\Pr((W_n,U_n)=(W,U))=\left\{
\begin{array}{cl}
\frac1{2|\mathcal{U}_n|} &
\textrm{if $U$, $W$ are both in $C$ with $|C\cap\mathcal{W}_n|=2$,}\\
\frac1{|\mathcal{U}_n|} &
\textrm{if $U$, $W$ are both in $C$ with $|C\cap\mathcal{W}_n|=1$,}\\
0 & \textrm{if $U\in C_1$ and $W\in C_2$ with $C_1\neq C_2$.}
\end{array}
\right.
\]
In view of the first part of Theorem~\ref{thm:combi}, the marginals
are uniformly distributed. The second part of Theorem~\ref{thm:combi}
gives the two inequalities.
\end{proof}

\subsection{Blossom trees and the conjugation of trees}
Theorem~\ref{thm:combi} is clearly analogous to Lemma~\ref{lem:cycle}
and~\ref{lem:height}. However we were not able to define directly
conjugacy classes on embedded trees. Instead we first construct an
encoding of embedded trees in terms of another family, blossom trees,
and then define conjugacy classes of blossom trees and prove
Theorem~\ref{thm:combi}.

\begin{figure}

\begin{center}
\includegraphics[width=9cm]{\rep 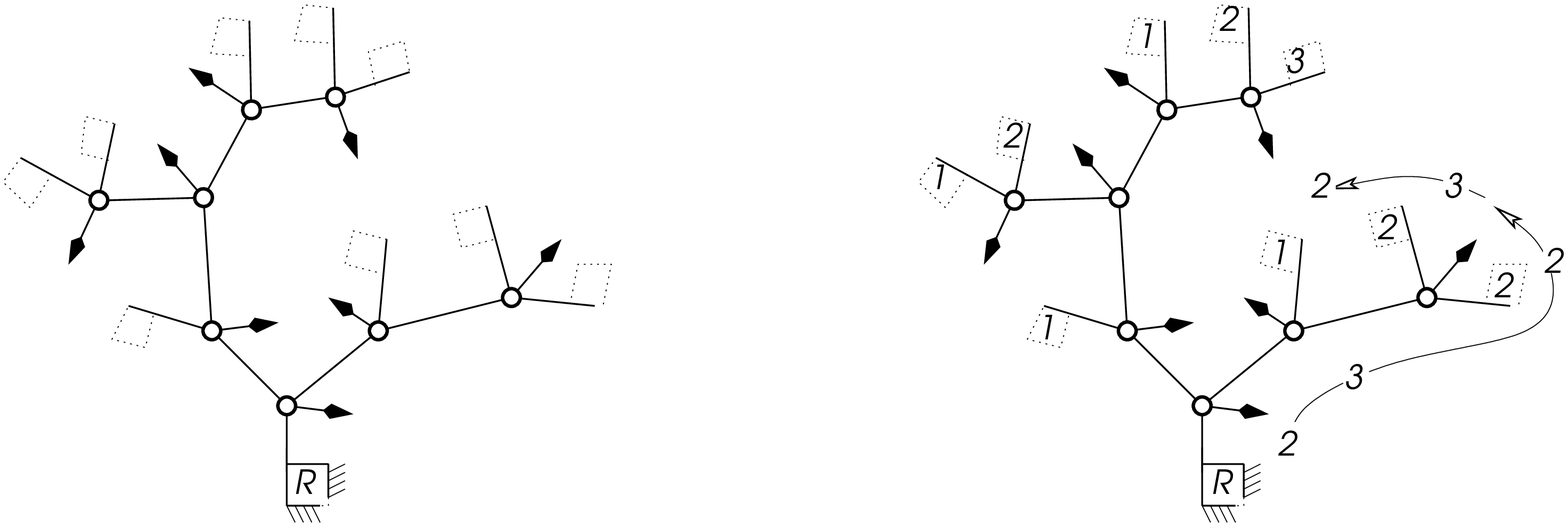}
\end{center}
\vspace{-5mm}
\caption{An example of blossom tree and its labelling.\label{fig:exofbin}}
\bigskip

\begin{center}
\includegraphics[width=11cm]{\rep 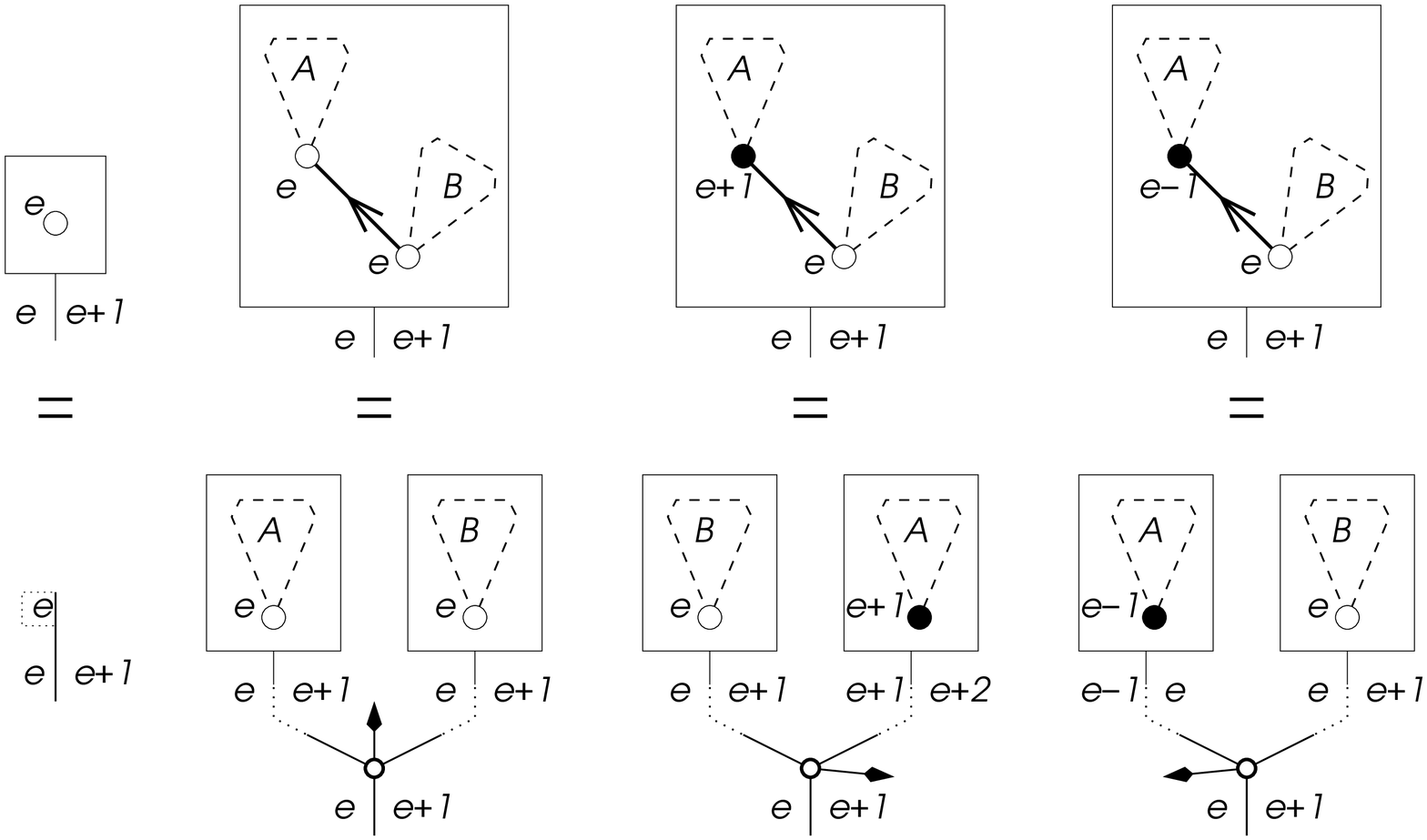}
\end{center}
\vspace{-5mm}
\caption{From embedded trees to blossom trees:
rules.\label{fig:toblossom}}
\bigskip

\begin{center}
\includegraphics[width=11cm]{\rep 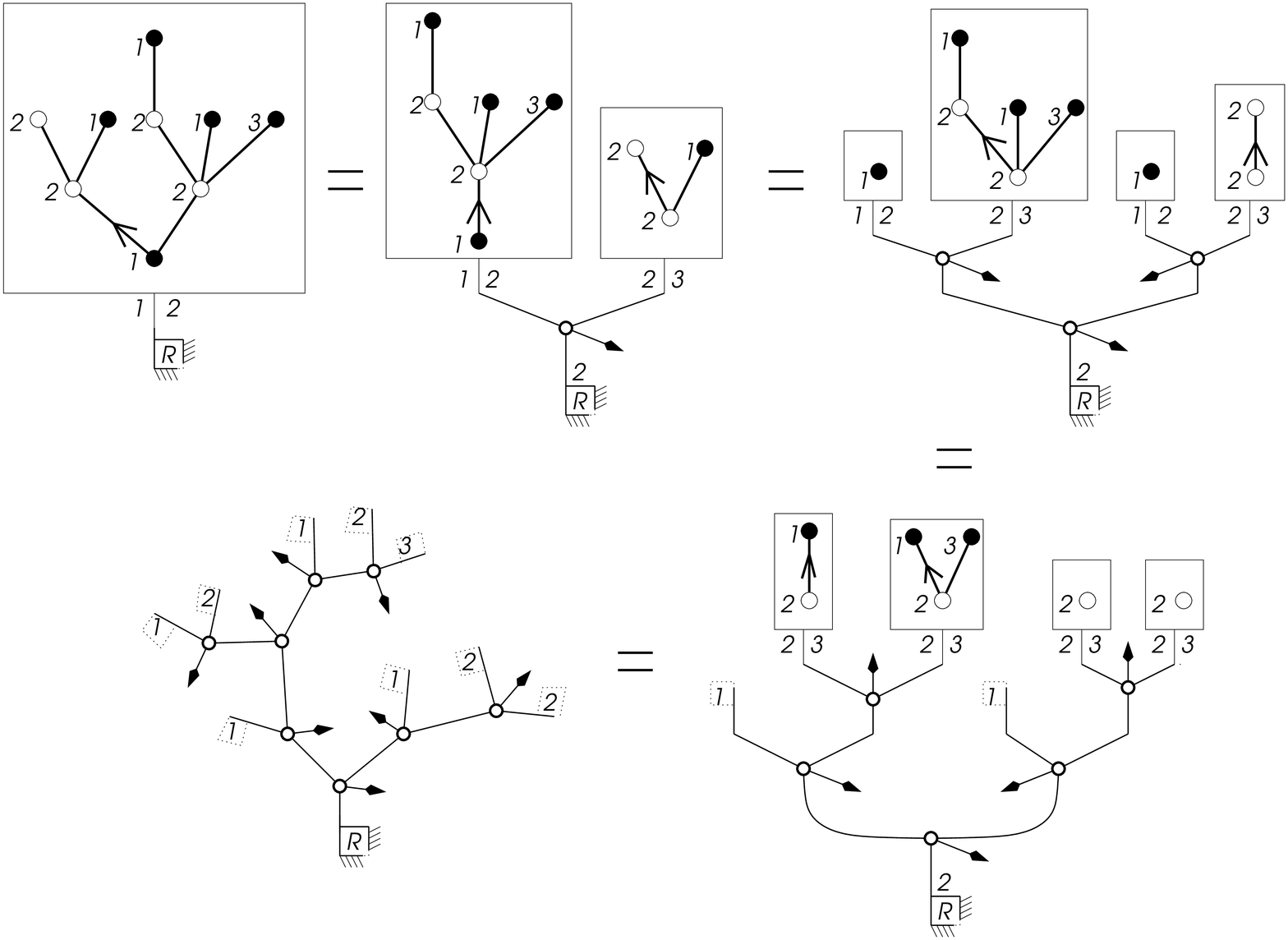}
\end{center}
\vspace{-5mm}
\caption{From embedded trees to blossom trees: an
example.\label{fig:example}}
\end{figure}

Following \cite{these-schaeffer}, let a \emph{blossom trees} be a
plane tree with the following properties:
\begin{itemize}
\item Vertices of degree one are of two types: arrows and flags. The
root of the a blossom tree is a flag, which is said to be
\emph{special}, as opposed to the other \emph{normal} flags.
\item All inner nodes have degree four and each of them is adjacent to
exactly one arrow.
\end{itemize}
Let $\mathcal{B}_n$ be the set of blossom trees with $n$ inner
nodes. By construction these trees have $n$ vertices of degree four
and thus $2n+2$ of degree one, that is $n$ arrows and $n+2$ flags.
The labelling of a blossom tree is given by the following
\emph{labelling process}:
\begin{itemize}
\item Start with \emph{current label} $2$ just after the root.
\item Turn around the border of the tree in counterclockwise direction.
\begin{itemize}
\item Each time an arrow is reached, the \emph{current label} is increased by
$1$.
\item Each time a flag is reached, the \emph{current label} is decreased by
one and then written on the flag.
\end{itemize}
\item Stop when the root flag is reached again (no label is written there).
\end{itemize}
This definition is illustrated by Figure~\ref{fig:exofbin}.

\begin{lem}
Embedded trees with $n$ edges are in one-to-one correspondence with
blossom trees with $n$ inner nodes with the same label distribution.
\end{lem}
\begin{proof}
In order to prove this lemma we work on a set of decorated blossom
trees: in these trees, the root flag is special and any flag with
label $e$ (as given by the labelling process) can either be empty or
be decorated by an embedded tree with root label $e$ (for $e\neq1$ the
immediate generalisation of embedded trees is meant). The
\emph{combined label distribution} of a decorated blossom tree counts
labels of decorations (embedded trees on flags) and of empty flags.  Examples
of decorated blossom trees are given in Figure~\ref{fig:example} (in
these figures, labels along edges indicate values taken by the current
label during the labelling process).

The first step of the encoding of an embedded tree consists in
writing it on the normal flag of the unique blossom tree with two
flags and no inner node (Figure~\ref{fig:example} top-left example).
Then the encoding is performed by recursively transforming the
decorated blossom tree according to the local rules of
Figure~\ref{fig:toblossom}. Each time the leftmost rule is applied to
a flag decorated by an embedded tree reduced to a vertex with label
$e$, this vertex is suppressed the flag becomes empty, with label $e$
(by construction, $e$ agrees with the labelling process; the combined
label distribution is left unchanged). When one of the other three
rules is applied, a new inner node is created while an edge of
embedded tree is suppressed. The relation between the position of the
created arrow and the root labels of the embedded trees grants that the
compatibility with the labelling process is preserved (observe that in
the middle rule subtrees have been switched for this purpose).

As long as there are decorated flags a rule can be applied. Once there
is no more decorated flag, a blossom tree is obtained. Rules are local
so that rules applied in distinct subtrees commute. As a
consequence the final blossom tree does not depend on the order in which
rules are applied.  Each rule is uniquely reversible so that the
encoding is bijective.
\end{proof}

\begin{figure}
\begin{center}
\includegraphics[width=12cm]{\rep 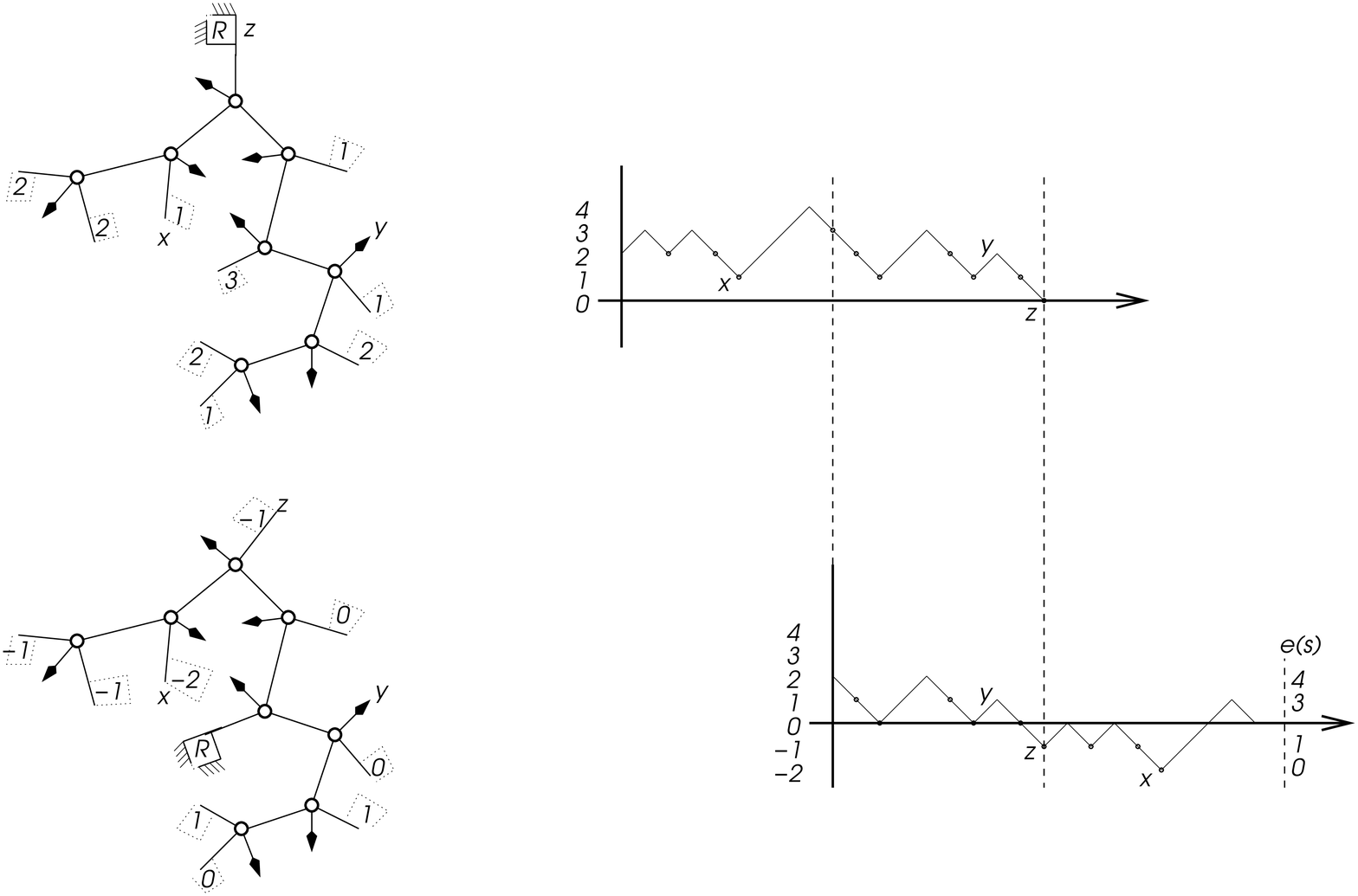}
\end{center}
\caption{The labelling process for two conjugated blossom trees. Three
steps $x$, $y$ and $z$ have been distinguished to illustrate the
correspondence: up step = arrow, down step = flag.
\label{fig:process}}
\end{figure}

\begin{proof}[Proof of Theorem~\ref{thm:combi}.]
The partition $\mathcal{U}_n\equiv\mathcal{B}_n=
\bigcup_{C\in\mathcal{C}_n}C$ is the partition of blossom trees in
conjugacy classes: two blossom trees $A$ and $B$ are in the same
\emph{conjugacy class} $C$ if $B$ is obtained from $A$ by first replacing
the root flag of $A$ by a normal flag and then choosing a new special
flag. This operation is called a \emph{cyclic shift} of the
tree.  In other terms each conjugacy class $C$ is the set of blossom
trees that can be obtained from a specific \emph{unrooted} blossom
tree (the flags of which are all normal) by selecting a special flag
(root flag) in all the possible ways.

Given a blossom tree $B$ with $n$ arrows and $n+2$ flags, the
evolution of the current label, while performing its labelling
process, is a walk $w_B$ with $n$ increments $+1$ and $n+2$ negative
increments $-1$, that starts from $2$, and whose last step, when the
process reaches again the root flag, is a negative increment. Upon
decreasing all labels by two, the walk $w_B$ is thus a walk of
$\mathcal{B}_{n,2}$ as defined in Section~\ref{sec:cycle}.

Moreover each cyclic shift of the tree $B$ is equivalent to the
corresponding cyclic shift of the walk $w_B$. Finally a blossom tree
encodes a well labelled tree if and only if all its labels are
positive, that is, if and only if the walk $w_B$ belongs to
$\mathcal{D}_{n,2}$ (upon decreasing all label by two).
The first statement of Theorem~\ref{thm:combi} is thus exactly the
cycle lemma (Lemma~\ref{lem:cycle}).

Finally, let us consider label distributions. Given $W$ a well
labelled tree and $U$ an embedded tree in the same conjugacy class of
trees, the corresponding walks $w_W$ and $w_U$ belong to the same
conjugacy class of walks. But the cumulated label distributions
satisfies $\widehat\lambda^{[W]}_k=\hat h_k(w_W)$ and
$\widehat\Lambda^{[U]}_k=\hat h_k(w_U)$ so that Lemma~\ref{lem:height}
gives the second statement of Theorem~\ref{thm:combi}.
\end{proof}

\section{Quadrangulations, Brownian snake and ISE}
\label{sec:proba}

\subsection{Encoding embedded trees by pairs of contour walks}
Let $\bar{\mathcal{U}}_n$ be the set of embedded trees with root label
zero instead of one. These trees, that are simply obtained from trees
of $\mathcal{U}_n$ by {shifting all labels down by one}, will be more
convenient for our purpose.

\begin{figure}
\begin{center}
\includegraphics[width=12cm]{\rep 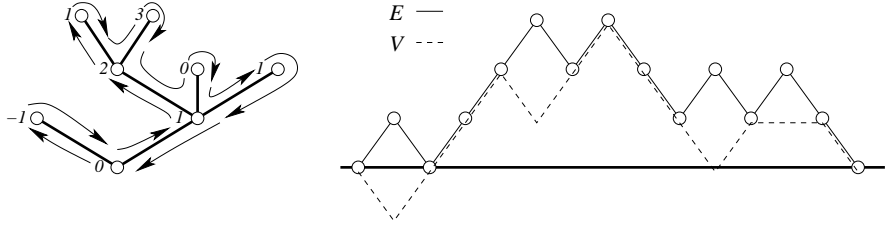}
\end{center}
\caption{Contour traversal and contour
pair $(E,V)$ of a tree.
\label{fig:traversal}}
\end{figure}

Let $U$ be an embedded tree of $\bar{\mathcal{U}}_n$ and consider the
following traversal of $U$, where traversing an edge takes unit time:
\begin{itemize}
\item At time $t=0$, the traversal arrives at the root.
\item If the traversal reaches at time $t$ a vertex $v_t$ having $k$ sons
for the $\ell$th time with $\ell\leq k$, its next step is toward the
$\ell$th son of $v_t$.
\item If the traversal reaches at time $t$ a vertex $v_t$ having $k$ sons
for the $(k+1)$th time, its next step is back toward the father of
$v_t$.
\end{itemize}
This traversal is called the \emph{contour traversal} because, as
exemplified by Figure~\ref{fig:traversal}, it turns around the
tree. In particular every edge is traversed twice (first away from and
then toward the root) and the complete traversal takes $2n$ steps.
The \emph{contour pair} of $U$ is then defined by the height
(\emph{i.e.}  distance to the root in the abstract tree), $E^{[U]}(t)$ and
label $V^{[U]}(t)$ of vertex $v_t$ traversed at time $t=0,\ldots,2n$.
(The path $E$ is often called the \textit{Dyck path} associated to the
tree $U$ \cite[Ch. 5]{Stanley}, or the \textit{contour process} in
\cite[Ch. I.3]{LEG}.)

The following proposition is immediate from the definition of contour pairs.
\begin{pro}\label{pro:contours}
The contour pair construction is a one-to-one correspondence between
$\bar{\mathcal{U}}_n$ (or $\mathcal{U}_n$) and the set
$\mathcal{EV}_{2n}$ of pairs of walks of length $2n$ such that:
\begin{itemize}
\item the walk $E$ is an excursion with increment $\pm1$ or Dyck path,
that is $E(0)=E(2n)=0$, $|E(t)-E(t+1)|=1$ and $E(t)\geq0$ for all
$t=0,\ldots,2n-1$;
\item the walk $V$ is a bridge with increment $\{-1,0,1\}$ or bilatere
Motzkin path, that is $V(0)=V(2n)=0$ and $(V(t)-V(t+1))\in\{-1,0,1\}$
for all $t$;
\item and the
consistency condition hold:\\
$\left(\; E(t)=E(t') \textrm{ and }
E(s)\geq E(t) \textrm{ for all } t<s<t'\;\right)
\;\Rightarrow\;V(t)=V(t').
$
\end{itemize}
\end{pro}
The excursion $E$ alone determines a unique unlabelled
rooted plane tree, while the walk $V$ describes one of the $3^n$
labelling of the tree encoded by $E$.  Recall that for an embedded
tree $U$, $\kappa(\epsilon)\in\{-1,0,1\}$ denotes the variation along
edge $\epsilon$ when traversed away from the root. In
particular if $\epsilon$ is traversed for the first time between time
$t$ and $t+1$ and for again between $t'$ and $t'+1$, then
\[\kappa(\epsilon)=V(t+1)-V(t)=V(t')-V(t'+1).\]
This local condition is equivalent to the consistency condition
of Proposition~\ref{pro:contours}.

\subsection{Random trees as random contour pairs}
\label{sec:contourpair}

Endow now $\bar{\mathcal{U}}_n$ with the uniform distribution and let
$(E^{(n)},V^{(n)})\equiv(E^{[U_n]},V^{[U_n]})$ denote the contour
pair of the random tree $U_n$. According to
Proposition~\ref{pro:contours}, the random contour pair
$(E^{(n)},V^{(n)})$ is uniformly distributed on $\mathcal{EV}_{2n}$
and $E_n$ is uniformly distributed on $\mathcal{E}_{2n}$, the set of
Dyck walks of length $2n$. More precisely, for all
$(E,V)\in\mathcal{EV}_{2n}$,
\[
\Pr((E^{(n)},V^{(n)})=(E,V))\;=\;\frac{1}{\frac{3^n}{n+1}{2n\choose n}},
\qquad
\Pr(E^{(n)}=E)\;=\;\frac{1}{\frac{1}{n+1}{2n\choose n}}.
\]

In order to state convergence results, let us now defined scaled
version of these random walks: given a random tree $U_n$ and its
contour pair $(E^{(n)},V^{(n)})$, let
\begin{eqnarray*}
e^{(n)}&=&\left(\frac{E^{(n)}(\lfloor2ns\rfloor)}{\sqrt{2n}}
\right)_{0\leq s\leq1}
\quad\textrm{and}\qquad
\hat W^{(n)}\;\;=\;\;
\left(\frac{V^{(n)}(\lfloor2ns\rfloor)}{(8n/9)^{1/4}}
\right)_{0\leq s\leq1}.
\end{eqnarray*}
The random variables $e^{(n)}$ and $\hat W^{(n)}$ take their values in
the Skorohod space $D([0,1],\mathbb{R})$ of c\`adl\`ag real functions
(right continuous with left limits).

As was proved by Kaigh \cite{KAIGH}, the scaled version $e^{(n)}$ of
the contour process converges weakly to the normalised Brownian
excursion~$e$. Our aim is to state an analogous result for the
random variable
\begin{eqnarray*}
X^{(n)}&\equiv&\parth{e^{(n)},\hat W^{(n)}},
\end{eqnarray*}
that takes its value in the Skorohod space $D([0,1],\mathbb{R}^2)$.

\subsection{A Brownian snake}
Let $e$ be the normalised Brownian excursion and
\[W=\parth{W_s(t)}_{0\le s\le 1,\ 0\le t\le e(s)}\]
be the Brownian snake with
lifetime $e$, as studied previously in
\cite{AALD,BCHS,DZ,DS,LEG,SERLET}.

\begin{figure}
\begin{center}
\includegraphics[width=10cm]{\rep 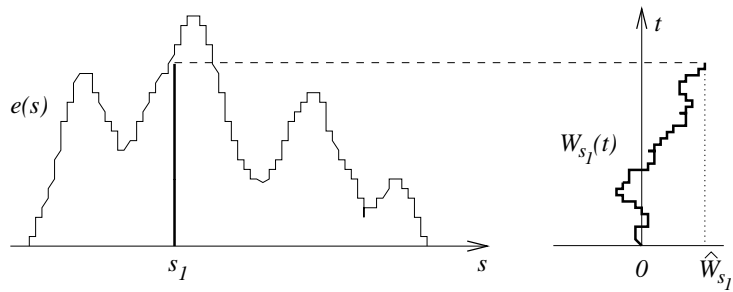}
\end{center}
\caption{Spacial extension of the snake at time $s_1$.
\label{fig:snake}}
\end{figure}

\begin{figure}
\begin{center}
\includegraphics[width=10cm]{\rep 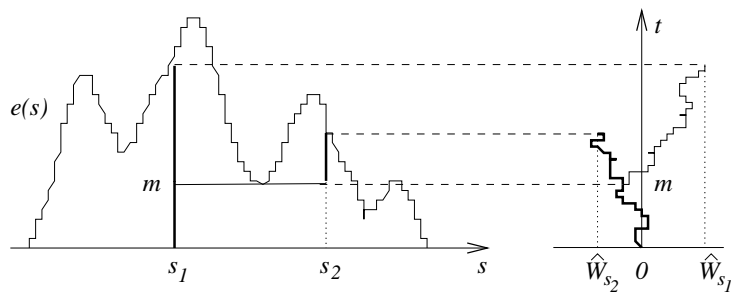}
\end{center}
\caption{Consistency of the snake between times $s_1$ and $s_2$.
\label{fig:snake-consist}}
\end{figure}

More precisely, the process $W$ can be defined as follows:
\begin{itemize}
\item
for all $0\leq s\leq 1$, $t\rightarrow W_s(t)$ is a standard Brownian
motion defined for $0\le t\le e(s)$ (see Figure~\ref{fig:snake});
\item the application $s\rightarrow W_s(.)$ is a
path-valued Markov process with transition function
satisfying:
     for $s_1<s_2$, and for $m=\inf_{s_1\le u\le s_2} \, e(u)$,
conditionally given $W_{s_1}(.)$ (see Figure~\ref{fig:snake-consist}),
\begin{itemize}
\item on the one hand we have that
\[
\parth{W_{s_1}(t)}_{0\le t\le m}=
\parth{W_{s_2}(t)}_{0\le t\le m},
\]
\item
and on the other hand $\parth{W_{s_2}(m+t)}_{0\le t\le
e(s_2)-m}$ is a standard Brownian motion starting from $W_{s_2}(m)$,
independent of $W_{s_1}(.)$.
\end{itemize}
\end{itemize}

The Brownian snake can be viewed as a branching Brownian motion, or as
an embedded continuum random tree (see \cite{AALD}).  More precisely
the excursion $e$ can be thought
of as the contour walk obtained by
contour traversal of a continuum random tree, while the snake
$W_s(\cdot)$ at times $s$ describes the embedding of the branch to the
root at time $s$.

\begin{figure}
\begin{center}
\includegraphics[width=10cm]{\rep 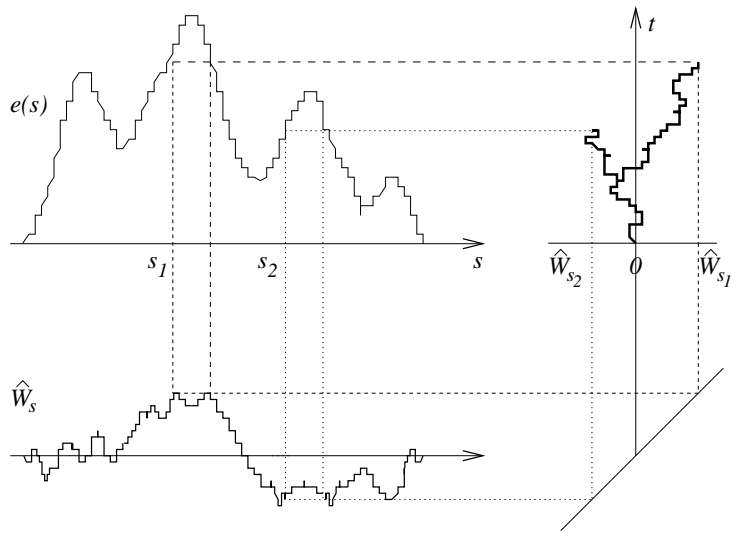}
\end{center}
\caption{The contour description $(e,\hat W_s)$: the excursion $e$
encodes the extension of the snake, the second walk
describes the horizontal position of its head.
\label{fig:contour-snake}}
\end{figure}

Instead of considering the full Brownian snake $W_s(t)$ we shall
concentrate, as we did in the discrete case, on its description by a
contour pair (or ``head of the snake'' description) $X=(X_s)_{0\leq
s\leq1}$, defined by  (see also Figure~\ref{fig:contour-snake})
\[
\hat W_s=W_s(e(s)), \quad X_s=\parth{e(s),\hat W_s}, \qquad
\textrm{for $0\leq s\leq1$.}
\]
In complete analogy with the discrete case, the full Brownian snake
can be reconstructed from its contour pair description since
$W_s(t)=\hat W_{\sigma(s,t)}$ where $\sigma(s,t)=\sup\{s'\leq s\mid
e(s')=t\}$. However we need only
and shall content with results in terms of $X$ (see \cite{MM} for a
complete discussion of the relation between the full snake and its
contour description).

\subsection{Integrated SuperBrownian Excursion}

Let $\mathcal{J}_n$ denote the empirical measure of labels of a random
embedded tree:
\[\mathcal{J}_n=\frac1n\sum_k\Lambda_k^{(n)}\ \delta_k.\]
Following Aldous \cite{AALD}, for any simple family of trees like our
embedded trees, $\mathcal{J}_n$ is expected to converge upon scaling
to a random mass distribution $\mathcal{J}$ supported by a random
interval $0\in[L,R]\subset\mathbb{R}$. This random measure $\mathcal
J$ is called Integrated SuperBrownian Excursion (ISE) by Aldous, in
view of its relation to $W$ through
\begin{equation}\label{bs2ise}
\int g\ d\mathcal{J}=\int_0^1\,g\parth{\hat W_s}\, ds,
\end{equation}
for any measurable test function $g$, see
\cite[Ch. IV.6]{LEG}.  In \cite{BCHS} the convergence of
$\mathcal{J}_n$ to $\mathcal{J}$ is proved for random embedded Cayley
trees. Although these trees are not exactly our random embedded
\emph{plane} trees, the proof could easily be adapted.

According to Corollary~\ref{cor:id} and Theorem~\ref{thm:coupling},
the radius $r_n$ is given by the width of the support of
$\mathcal{J}_n$.
However the weak convergence of $\mathcal{J}_n$ to $\mathcal{J}$,
as obtained in \cite{BCHS} is not sufficient for our purpose since
$r=R-L$, the width of the support of $\mathcal J$, is not a continuous
functional of the measure $\mathcal J$.

\subsection{Convergence of snakes} Instead of weak convergence of
$\mathcal J_n$
to $\mathcal J$, we
shall thus prove in Section~\ref{sec:conv} the following stronger
result.
\begin{theo}\label{thm:converge} The scaled contour pair
$X^{(n)}$  converges weakly  to $X$ in $D([0,1],\mathbb{R}^2)$.
\end{theo}
This theorem establishes weak convergence of the scaled contour (or
head of the snake) description of embedded trees to the head of the
snake description of the Brownian snake with lifetime $e$.  We
moreover obtain a deviation bound for the maximal extension of the
snake $\hat W^{(n)}_s$.
\begin{pro}\label{pro:devbound}
There exists $y_0>0$ such that for all $y>y_0$ and $n$,
\[
\mathbb{P}\left(\sup_{0\leq s\leq1}\hat
W^{(n)}_s\,>\,(8/9)^{1/4}y\right)\;\leq e^{-y}.
\]
\end{pro}

Theorem~\ref{thm:converge} was independently
obtained by Marckert and Mokkadem \cite{MM}. They extend the
convergence result to the explicit full description $(W_s(t))_{s,t}$
but their alternative proof does not provide the exponential bound
of Proposition~\ref{pro:devbound}.

\subsection{The radius of a random quadrangulation and the width of ISE}
According to Corollary \ref{cor:id} and to
Theorem \ref{thm:coupling}, the radius $r_n$ of the
           quadrangulation  corresponding  to $U_n$ satisfies
\[\left\vert
         \parth{8/9}^{1/4}\parth{\sup_{0\le s\le 1} \hat W^{(n)}_s
-\inf_{0\le s\le 1} \hat W^{(n)}_s}-\,n^{-1/4}\,r_n\right\vert\le
3n^{-1/4}.\]
Theorem~\ref{thm:converge} and
Proposition~\ref{pro:devbound} thus prove the
conjecture $\smash{\mathbb{E}(r_n)=\Theta(n^{1/4})}$ and lead to a
much more precise characterization:
\begin{cor}
\label{radcv}
The random variable $n^{-1/4}\ r_n$  converges weakly  to
$(8/9)^{1/4}\ r$, in which
          \[r=\sup_{0\le s\le 1}\hat W_s-\inf_{0\le s\le 1}\hat W_s.\]
Furthermore, convergence of all moments holds true.
\end{cor}
In view of Relation~(\ref{bs2ise}) the random variable $r$ is also the
width of ISE process~$\mathcal J$.

\subsection{The profile and a CRM}
Actually, Theorems~\ref{thm:welllab} and~\ref{thm:coupling}  suggest that
not only the scaled radius but the full scaled profile converges (at
least in distribution) to the ISE mass distribution. More precisely,
define the distribution function $F(x)$ of the translated ISE by
\[W_{\min}=\inf_{0\le s\le 1}\hat W_s,
\hspace{0,5cm}F(x)=\mathcal J\parth{(-\infty,W_{\min}+x]}
=\mathcal{J}\parth{[W_{\min},W_{\min}+x)]},\]
and the scaled distribution function of the profile of random
quadrangulations by
\[F_n(x)
=\frac 1{n+1}\  \widehat \lambda^{(n)}_{\lfloor(8n/9)^{1/4}x\rfloor}
=\frac 1{n+1}\  \widehat H^{(n)}_{\lfloor(8n/9)^{1/4}x\rfloor}.\]
where $\widehat\lambda^{(n)}_k$ is the cumulated distribution of
labels of a random well labelled tree (as defined in
Section~\ref{sec:welllab}) and $\widehat H^{(n)}_k$ is the cumulated
profile of a random quadrangulation (as defined in
Section~\ref{sec:def}).

Then we prove the following corollary of Theorems~\ref{thm:welllab},
\ref{thm:coupling}, \ref{thm:converge} and Corollary~\ref{radcv}.

\begin{cor}
\label{profile}
The scaled profile $F_n$  converges weakly  to
$F$ in $D([0,+\infty),\mathbb{R})$.
\end{cor}

A natural conjecture is that there is a continuum analogue to
Theorem~\ref{thm:welllab} that allows to define from ISE a Continuum
Random Map (CRM), such that the properties of scaled distances in
random quadrangulations (distances between arbitrary pairs of points,
not only with respect to a basepoint) would be described by the
properties of distance in the CRM. In view of the interpretation of
random quadrangulations as 2d Euclidean pure quantum geometries, this
CRM might be considered as a natural candidate model of continuum 2d
pure quantum geometry. We plan to discuss this connection further in a
subsequent paper.

\section{Random embedded trees and the Brownian snake}
\label{sec:conv}

In this section we prove Theorem~\ref{thm:converge}.  Finite
dimensional density functions are first calculated
(Section~\ref{sec:fidis}).  Section~\ref{sec:bound} then provides the
deviation bound for the maximum of the label walk. Finally tightness
is proved using the previous bound (Section~\ref{sec:tight}). The
theorem then follows from standard results on weak convergence in the
space $D([0,1],\mathbb{R}^2)$ \cite{BILL}.

\subsection{Finite dimensional density functions}
\label{sec:fidis}

        From now on in this section, $p$ and ${\tau}=(\tau_1,\dots,\tau_p)$ are
fixed with $0<\tau_1<\cdots<\tau_p<1$, and we prove the following
finite dimensional density convergence result.
\begin{pro}\label{pro:fidis}
The sequence of random variables
\[X^{(n)}(\tau)=
\left(\frac{E^{(n)}(\lfloor2\tau_in\rfloor)}{(2n)^{1/2}},\;
\frac{V^{(n)}(\lfloor2\tau_in\rfloor)}{(8n/9)^{1/4}}
\right)_{1\le i\le p}=\left(e^{(n)}(\tau_i),\;
\hat W^{(n)}_{\tau_i}\right)_{1\le i\le p}
\]
weakly converges to  $X(\tau)=\left(e(\tau_i),\;
\hat W_{\tau_i}\right)_{1\le i\le p}$,
that is,
\[\lim_n\mathbb{E}\left[\Phi\left(X^{(n)}(\tau)\right)\right]=
\mathbb{E}\left[\Phi\left(X(\tau)\right)\right]\]
or any bounded   continuous
function $\Phi$ on $\mathbb{R}^{2p}$.
\end{pro}

For this aim it will be convenient to prove first the weak convergence
of
\[Y^{(n)}(\tau)=\left(\left(e^{(n)}(\tau_i)\right)_{1\le i\le p},
\Big(
\inf_{[\tau_i,\tau_{i+1}]}e^{(n)}
\Big)_{1\le i\le p-1}\right),\]
then that of
\[Z^{(n)}(\tau)=\left(Y^{(n)}(\tau),\;\left(
\hat W^{(n)}_{\tau_i}\right)_{1\le i\le p},
\left(\mu^{(n)}_{i}\right)_{1\le i\le p-1}\right),\] in which
$\mu^{(n)}_i$ is defined by
\[
\mu^{(n)}_i\;=\;\hat W^{(n)}_{\tau'_i} \qquad \textrm{ for any }\quad
\tau'_i\in\mathop{\mathrm{arginf}}_{[\tau_i,\tau_{i+1}]}e^{(n)}.
\]
(The consistency condition of Proposition~\ref{pro:contours} grants
     that $\mu^{(n)}_i$ is indeed independent of the exact choice of
     $\tau'_i$ in
     $\mathop{\mathrm{arginf}}_{[\tau_i,\tau_{i+1}]}e^{(n)}$.)

The weak convergence of $Y^{(n)}(\tau)$ follows from \cite{KAIGH},
as a special case, but does not fill our needs.
In the next section we prove a local limit
      theorem for finite dimensional distributions of
$e^{(n)}$ which is not a consequence of \cite{KAIGH}.
Once the local limit theorem for $Y^{(n)}(\tau)$ is proved, we recall
its interpretation in terms of trees, using the key notion of
\emph{shape}. This leads us to split $\mathbb{R}^{4p-2}$ in $(p-1)!$
regions $(R_\sigma)_{\sigma\in\mathbb{S}_{p-1}}$ and to prove weak
convergence of $Z^{(n)}(\tau)$ separately on each region.  Finally we
identify the limit and the weak convergence of $X^{(n)}(\tau)$ follows
from that of $Z^{(n)}(\tau)$.

\subsubsection{A local limit theorem for  $Y^{(n)}(\tau)$.}

The characteristics of the walk $e^{(n)}$ we are interested in
        are contained in  the
sequence of successive heights at (resp.  minima between) the
$\tau_i$. Let $(x,m)=\parth{(x_i)_{1\le i\le p},(m_i)_{1\le i\le p-1}}$
be a typical value of $Y^{(n)}(\tau)$, that is,  $x_i\sqrt{2n}$
         and $m_i\sqrt{2n}$  are   integers, and they satisfy
\[
\inf(x_i,x_{i+1}) \geq m_i, \hspace{1cm}1\le i\le p-1.
\]
For the first result, we need a more handy
parametrisation, by the successive up and down
relative variations: let $\gamma_0=x_1$, $\beta_p=x_p$,
\begin{equation}\label{rel:bgxm}
\gamma_i=x_{i+1}-m_{i},\qquad \beta_i=x_i-m_i, \qquad
\textrm{for $i=1,\ldots,p-1$},
\end{equation}
and by convention $\beta_0=\gamma_p=0$.
\begin{pro}\label{pro:fdd-shape}
Let $K$ be a compact subset of
\[
\left\{(x,m)\in\mathbb R^p\times \mathbb R^{p-1}
\ \left\vert\ 0<m_i<\inf(x_i,x_{i+1}),1\le i\le p-1\right.\right\},
\]
that is, of the domain of definition of coherent values of the $x_i$ and $m_i$.
Then,
uniformly for $(x,m)\in K\cap \big((2n)^{-1/2}\mathbb{Z}\big)^{2p-1}$,
\begin{eqnarray*}
\mathbb{P}\left(Y^{(n)}(\tau)=\left(
x,m\right)\right)
\,=\,(2n)^{-2\delta}\cdot \zeta(x,m)\cdot
\parth{1+\mathcal{O}_K(n^{-1/2})},
\end{eqnarray*}
where $\delta=(2p-1)/4$, and $\zeta(x,m)$ reads, in terms of
the $\beta_i$ and $\gamma_i$ as given by relations~(\ref{rel:bgxm}):
\begin{eqnarray*}
\zeta(x,m)&=&2^{2p}\,(2\pi)^{-p/2}
\prod_{i=0}^p\frac{(\beta_i+\gamma_i)
e^{-\frac{(\beta_i+\gamma_i)^2}{2(\tau_{i+1}-\tau_{i})}}}
{(\tau_{i+1}-\tau_{i})^{3/2}}.
\end{eqnarray*}
\end{pro}
The notation $\mathcal{O}_K$ is used to stress the fact that the error term
is uniform for fixed $K$.
\begin{proof}
Let, for all $n$ and $a$ non negative integers,
\begin{eqnarray}
\nonumber C(n;a)&=&\frac{a+1}{n+1}{n+1\choose (n-a)/2}\\
\label{rel:asympt}&=&\frac{2^{n+1}}{\sqrt{2\pi} n}\
\frac{(a+1)}{\sqrt{n}}\,e^{-\frac{a^2}{2n}}
\left(1+O({\textstyle\frac{a}n})
\right),
\end{eqnarray}
where the error term is uniform for $a=O(n^{1/2})$.

The Catalan numbers $C(2n;0)$ are well known to give the cardinality of
$\mathcal{E}_{2n}$. More generally, the reflection principle proves
that $C(n;a)$ is the number of meanders with increments $\pm1$ (aka
left factors of Dyck walks), that have length $n$ and end at height
$a$.  Given non negative numbers $b$ and $c$, $C(n; b+c)$ is also the
number of walks with increments $\pm1$ that have length $n$, minimum
$-b$ and final height $c-b$, as follows from a decomposition at first
and last passage at the minimum.

We thus have
\[
\mathbb{P}\left(Y^{(n)}(\tau)=(x,m)\right)\;=\;
\frac{\prod_{i=0}^pC\left(\lfloor2n\tau_{i+1}\rfloor-
\lfloor2n\tau_{i}\rfloor\,;\;(2n)^{1/2}
\beta_i+(2n)^{1/2}\gamma_i\right)}{C(2n;0)}.
\]
Combined with (\ref{rel:asympt}), it  yields
Proposition~\ref{pro:fdd-shape}.
\end{proof}
The expression
\begin{equation}\label{def:zeta}
\zeta(x,m)\;=\;2^{2p}(2\pi)^{-p/2}\
\prod_{i=0}^p\ \ \frac{(\beta_i+\gamma_i)
e^{-\frac{(\beta_i+\gamma_i)^2}
{2(\tau_{i+1}-\tau_i)}}}{(\tau_{i+1}-\tau_{i})^{3/2}}\
\mathbf{1}_{\beta_i\ge0}
\mathbf{1}_{\gamma_i\ge0}
\mathbf{1}_{m_i\ge0}
\end{equation}
where the $\beta_i$ and $\gamma_i$ are given by
relations~(\ref{rel:bgxm}), is the expected limit density of
probability for the $x_i$ and $m_i$: in particular it is
coherent with the density of evaluations of the normalized Brownian
excursion at $p$ points, and of the $p-1$ minima between them, as
given in \cite{SERLET}.

\subsubsection{Shapes}\label{subs:shapes} Let us now define a
\emph{shape} to be a rooted plane tree $T$ with $q$ edges, that we
call \emph{superedges} to distinguish them from edges of
embedded trees.  We shall endow each superedge with a length: let
$\eta_1$, \ldots, $\eta_q$ denote the edges of $T$ in \textit{prefix order}
(i.e. the order induced by first visits in contour traversal) and
$\mathcal{L}(\eta_i)$ the length of $\eta_i$.

\begin{figure}
\begin{center}
\includegraphics[width=13cm]{\rep 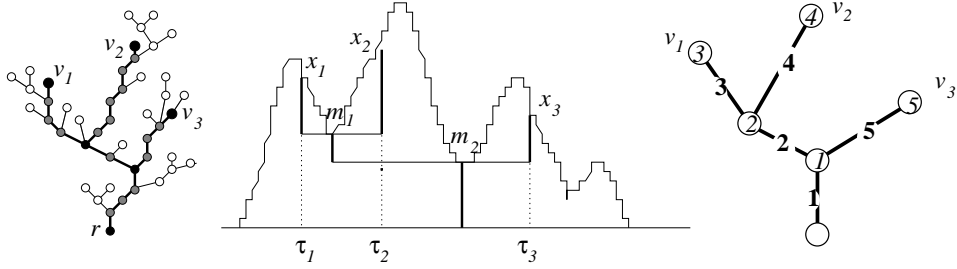}
\end{center}
\caption{A random tree $U_n$, the decomposition of its normalised
contour $e^{(n)}$ at times $(\tau_1,\tau_2,\tau_3)$ and at minima in
between, and the extracted shape (prefix order shown on
edges and vertices). In particular: $\mathcal{L}(\eta_4)=7=(x_2-m_1)\sqrt{2n}$.
\label{fig:shape}}
\end{figure}

Given the $p$ normalised times $0<\tau_1<\cdots<\tau_p<1$ and $U$ an
embedded tree of $\bar{\mathcal{U}}_n$, let us extract a shape $T$
and, for each superedge $\eta$ of $T$, a length
$\mathcal{L}(\eta)$. In order to do this, let
$t_i=\lfloor2n\tau_i\rfloor$ and consider $v_{t_i}$ the vertex visited
at time $t_i$ by the contour traversal. The $p$ \emph{fixed vertices}
$v_{t_i}$ and the root $r$ span a minimal subtree of $U$ (the union of
the branches from $v_{t_i}$ to $r$, see Figure~\ref{fig:shape}, left
hand side). Apart from the $v_{t_i}$, this subtree contains vertices
of two types: \emph{branchpoints} that have at least two sons in the
subtree (black vertices in Figure~\ref{fig:shape}), and \emph{smooth
vertices} that have exactly one son in the subtree (grey
vertices). Let us define the shape $T$ by taking the fixed vertices,
the root and the branchpoints as vertex set, and the paths connecting
these vertices as superedges (Figure~\ref{fig:shape}, right hand
side). A superedge $\eta$ of $T$ is thus by definition a set of edges
of $U$ and we let $\mathcal{L}(\eta)=|\eta|$, which is just the length
of the path $\eta$ in the tree $U$ (for instance, superedge $\eta_4$
in Figure~\ref{fig:shape} is made of $7$ edges).

Observe now that the shape and superedge lengths extracted from a
random tree $U_n$ completely determine $Y^{(n)}(\tau)$.  To check this
assertion, let us assume that the extraction yields a shape $T$ and
superedges $\eta_i$ of normalised length $\ell(\eta_i)$ with
$\mathcal{L}(\eta_i)=\ell(\eta_i)\sqrt{2n}$, $i=1,\dots, q$.  Let
$A_i$ be the set of superedges on the unique path from $v_{t_i}$ to
the root, so that $\sum_{\eta\in A_i}\ell(\eta)=e^{(n)}(\tau_i) =
x_i$.  The contour traversal of the shape $T$ starts from the root
$v_{t_0}=v_{t_{p+1}}$ and reaches successively each leaf $v_{t_1}$,
\dots, $v_{t_p}$. From $v_{t_i}$ to $v_{t_{i+1}}$ a set $B_i$ of
superedges are first traversed toward the root, followed by a set
$C_i$ of superedges that are traversed away from the root (so that
$B_i=A_i\setminus A_{i+1}$ while $C_i=A_{i+1}\setminus A_i$).  Then
$\smash{\sum_{\eta\in B_i}\ell(\eta)}=\beta_i$ and
$\smash{\sum_{\eta\in C_i}\ell(\eta)}=\gamma_i$ are the lengths of
these journeys, with $\beta_i$ and $\gamma_i$ given by
Relations~(\ref{rel:bgxm}).  In particular the normalised lengths
$\ell(\eta)$ of superedges are all of the form
\begin{equation*}
(x_i-m_i),\hspace{7mm}(x_{i+1}-m_i)\hspace{7mm}\textrm{or}\hspace{1cm}
\vert m_j-m_i\vert \ \ \textrm{with} \ \ j< i.
\end{equation*}
In the latter case, $m_j$ has to be a record, that is for $j<k<i$,
$m_k>\sup(m_i,m_j)$.
These relations are exemplified by Figure~\ref{fig:shape}.

Conversely, as explained in \cite[Ch. 3]{LEG} or
\cite[Section~2]{SERLET}, the $x_i$ and $m_i$
(or the
$\beta_i$ and $\gamma_i$ as defined by Relations~(\ref{rel:bgxm}))
exactly determine the shape and superedge lengths of $U_n$.

\begin{figure}
\begin{minipage}{0.45\textwidth}
\includegraphics[width=5.5cm]{\rep 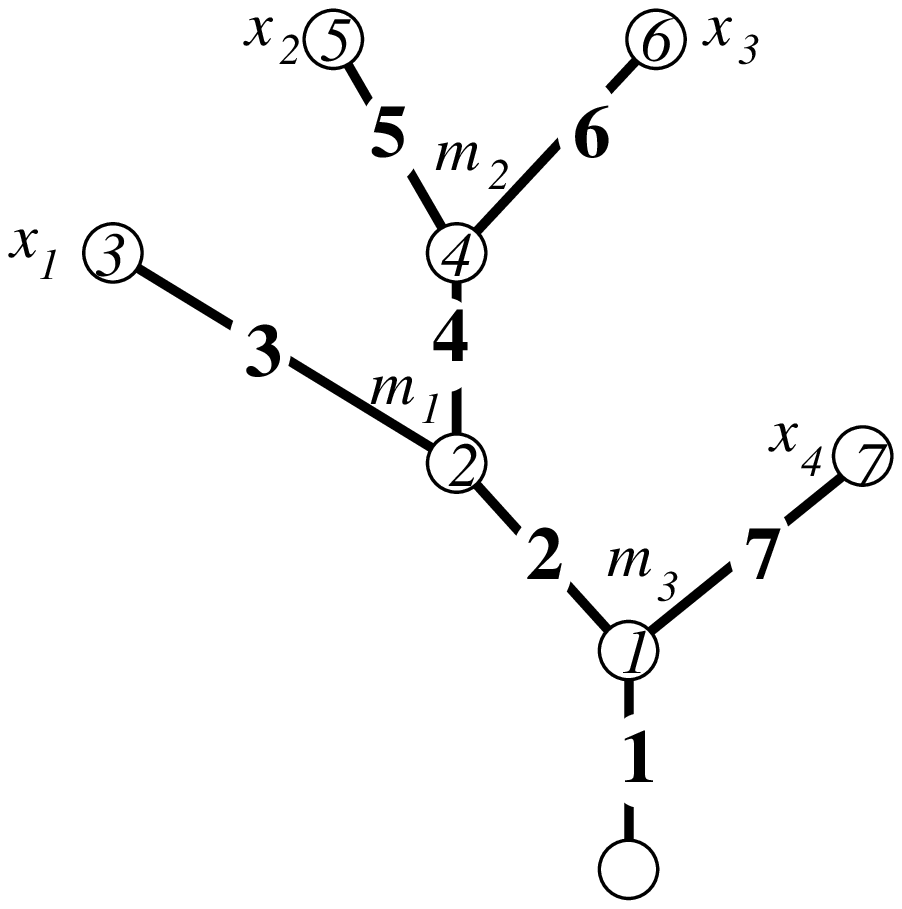}
\end{minipage}
\hspace{10mm}
\begin{minipage}{0.45\textwidth}
{\Large $
\begin{array}{c|ccccccc}
&\eta_1\!&\!\eta_2\!&\!\eta_3\!&\!\eta_4\!&\!\eta_5\!&
\!\eta_6\!&\!\eta_7\\\hline
m_3&  1&0&0&0&0&0&0\\
m_1&  1&1&0&0&0&0&0\\
x_1&  1&1&1&0&0&0&0\\
m_2&  1&1&0&1&0&0&0\\
x_2&  1&1&0&1&1&0&0\\
x_3&  1&1&0&1&0&1&0\\
x_4&  1&0&0&0&0&0&1
\end{array}
$
}
\end{minipage}
\caption{A shape $T$ with prefix ordering of superedges and vertices,   and
the matrix $M_T$ mapping superedge lengths onto $(x,m)$.
\label{fig:matrix}}
\end{figure}
>From the limit density $\zeta(x,m)$ of the previous section, the
probability that $m_i=m_j$ for some $i\neq j$ is seen to tend to zero as $n$
goes to infinity. This implies that, with probability tending to one as
$n$ goes to infinity, the shape is a binary tree and $q=2p-1$ (indeed
the existence of a branchpoint of larger degree corresponds to the
equality of two minima $m_i$ and $m_j$).  From now on we thus restrict
our attention to binary shapes $T$.

Associated to any binary shape $T$, there is a matrix $M_T$ with size
$2p-1\times 2p-1$ and entries zero or one, that sends the $2p-1$
normalized lengths of superedges $(\ell(\eta_i))$ on the $(x_i,m_i)$s.
Moreover, provided the lengths of superedges on the one hand, and the
$x_i$ and $m_i$ on the other hand, are sorted according to the prefix
order, the matrix $M_T$ is lower triangular with ones on the
diagonal. As a consequence, $M_T$ is \textit{Lebesgue measure
preserving}. Indeed the $k$th vertex in the prefix order is reached by the
$k$th edge and all other edges on the path to the root have already
been visited.

Finally let us consider labels of $U_n$: the variation of
label along edges extends to superedges, upon setting for any
syperedge $\eta$
\[
\kappa(\eta)=\sum_{\epsilon\in\eta}\kappa(\epsilon).
\]
Let $\tilde M_T$ be the double action of $M_T$ on $\mathbb R^{4p-2}$.
Since the matrix $M_T$ also sends the normalized increments
$(k_i)_{1\leq i\leq 2p-1}$ of labels along the $2p-1$ superedges onto
the normalized labels $(w_i,\mu_i)$ of the $2p-1$ vertices of the
tree,  $\tilde M_T$ describes the one-to-one correspondance between
shapes equiped with superedges' lengths and label
variations on the one hand, and $Z^{(n)}(\tau)$ on the other hand.
More precisely, $\tilde M_T$ describes  the restriction of
this correspondance to shape $T$.

\subsubsection{A local limit theorem for $Z^{(n)}(\tau)$ in a fixed region.}
\label{subs:Zfidis}

As already observed, the limit law
of $Y^{(n)}(\tau)$ charges only the region
\[R=\left\{(x,m,w,\mu)\in(\mathbb R^p\times \mathbb R^{p-1})^2
\ \left\vert\ 0<m_i<\inf(x_i,x_{i+1}),1\le i\le p-1\right.\right\}.\]
For a given permutation $\sigma$ on $p-1$ symbols, define
\[R_{\sigma}=\left\{\left.(x,m,w,\mu)\in R
\ \right\vert\ 0<m_{\sigma(1)}<m_{\sigma(2)}
<\dots<m_{\sigma(p-1)}\right\}.\]
Observe that the shape is constant on $R_{\sigma}$ and denote
it $T_\sigma$. (Conversely a
given shape may appear in many different regions
$R_{\sigma}$, as there are only $C(p-1)$ different shapes.) Let
\[R_{\sigma,\varepsilon}=\left\{\left.(x,m,w,\mu)\in R_{\sigma}
\ \right\vert\
d\left((x,m,w,\mu),\partial  R_{\sigma}\right)
\ge\varepsilon \right\}.\]
(The distance is the usual distance in $\mathbb{R}^{4p-2}$ and
$\partial R_\sigma$ denote the boundary of $R_\sigma$, which is clearly
a finite union of closed $(4p-3)$-dimensional cones.)

Let $\Delta(n,K)$ denote
the set of possible values of $Z^{(n)}(\tau)$ that
belong to a compact subset $K$ of $\mathbb R^{4p-2}$,
that is
\[
\Delta(n,K)\;=\;K\cap\parth{\parth{(2n)^{-1/2}\mathbb Z}^{2p-1}
\times\parth{\parth{8n/9}^{-1/4}\mathbb Z}^{2p-1}}.
\]
Furthermore, let  $\xi$ (resp. $f$) be defined,
on $\mathbb R^{2p-1}\times \mathbb R^{2p-1}$
(resp. $(\mathbb R^p\times \mathbb R^{p-1})^2$), by
\begin{eqnarray}\label{equ:defxi}
\xi(\ell,k)&=&(2\pi)^{-2\delta}
\prod_{i=1}^{2p-1}\;\ell_i^{-1/2}\;
\exp\left(-\frac{k_i^2}{2\,\ell_i}\right),\\
\nonumber
\xi_T&=&\xi\circ \tilde M_T^{-1},\\
\nonumber f(x,m,w,\mu)&=&\zeta(x,m)\ \sum_{\sigma}\xi_{T_\sigma}(x,m,w,\mu)
\;\cdot\; \mathbf{1}_{R_{\sigma}}(x,m).
\end{eqnarray}
The function $f$ is exactly the density of $Z(\tau)$, that is, the
density of the
evaluation of $X$ at the $p$ points $\tau_i$ and at the $p-1$ minima
of $e$ between them. This density was described in \cite[Propositions
2 \& 3]{SERLET}. The function $\xi_{T_\sigma}$  is the conditional density
of the labels $(w,\mu)$ given  $(x,m)$.
We shall prove the following local limit law for $Z^{(n)}(\tau)$.
\begin{lem}\label{lem:locZ}
Let $K$ be a compact subset of
$R_{\sigma,\varepsilon}$, $\epsilon>0$.  Then, uniformly for
$(x,m,w,\mu)\in \Delta(n,K)$,
\begin{eqnarray*}
\lefteqn{\mathbb{P}\left(Z^{(n)}(\tau)=(x,m,w,\mu)\right)}\\
&=&
(2n)^{-2\delta}(8n/9)^{-\delta}
\ \ f(x,m,w,\mu)\ \ \parth{1+\mathcal
O_K\parth{n^{-1/2}}},
\end{eqnarray*}
\end{lem}
\begin{cor}\label{cor:weakZ}
Let $K$ denote a compact subset of
      $R_{\sigma,\varepsilon}$.
For any uniformly continuous function $\Phi$ with support
      $K$,
\begin{eqnarray*}
\lim_n\mathbb{E}\left[\Phi\left(Z^{(n)}(\tau)\right)\right]&=&
\mathbb{E}\left[\Phi\left(Z(\tau)\right)\right].
\end{eqnarray*}
\end{cor}
\begin{proof}[Proof of Lemma~\ref{lem:locZ}]
Our aim is to compute for $(x,m,w,\mu)\in \Delta(n,K)$ the probability
\begin{eqnarray*}
\lefteqn{\mathbb{P}\left(Z^{(n)}(\tau)=(x,m,w,\mu)\right)}\\
&=&
\mathbb{P}\left(Y^{(n)}(\tau)=(x,m)\right)\cdot
\mathbb{P}\left(Z^{(n)}(\tau)=(x,m,w,\mu)
\mid Y^{(n)}=(x,m)\right).
\end{eqnarray*}
>From Proposition~\ref{pro:fdd-shape}, we already have
\begin{eqnarray*}
\mathbb{P}\left(Y^{(n)}(\tau)=(x,m)\right)
&=&
(2n)^{-2\delta}\ \ \zeta(x,m)\ \parth{1+O\parth{n^{-1/2}}}.
\end{eqnarray*}
As discussed in the previous section, the normalised
      lengths $\ell_i=(2n)^{-1/2}\mathcal{L}_i$ of
superedges, obtained from $(x,m)$ through $M_T$, are of the form
\begin{equation*}
(x_i-m_i),\hspace{7mm}(x_{i+1}-m_i)\hspace{7mm}\textrm{or}\hspace{1cm}
\vert m_j-m_i\vert \ \ \textrm{with} \ \ j< i.
\end{equation*}
In particular the fact
that $(x,m,\cdot,\cdot)\in K\subset
R_{\sigma,\varepsilon}$ for some
$\epsilon>0$ grants that these normalized lengths are uniformly
bounded away from~0.  In turn the variation $\kappa^{(n)}(\eta)$ of
labels along any superedge $\eta$ is the sum of at least $\mathcal
O_{\varepsilon}(\sqrt{n})$ i.i.d. uniform random variables on
$\{1,0,-1\}$. Therefore we can apply uniform bounds for the local
limit theorem \cite[pages 189--197]{PET}: if $S_n$ denotes the sum of
$n$ i.i.d.  random variables uniform on $\{+1,0,-1\}$, we have
\[\mathbb{P}\left(S_n=\kappa\right)
=\ \sqrt{\frac{3}{4\pi n}}\ e^{-3\kappa^2/4n}+O\left(\frac 1n\right).
\]
This allows us to calculate the probability that the variations of
labels along superedges $\kappa^{(n)}(\eta_i)$ are equal to
$\kappa_i=k_i\parth{8n/9}^{1/4}$, $i=1,\ldots,2p-1$:
uniformly for $(x,m,w,\mu)\in\Delta(n,K)$,
\begin{eqnarray*}
\lefteqn{\mathbb{P}\left(\kappa^{(n)}(\eta_i)=\kappa_i,\;{1\leq i\leq
2p-1}\;\Big\vert\; Y^{(n)}(\tau)=(x,m) \right)}\\
&=&\ \prod_{i=1}^{2p-1}\left(\sqrt{\frac{3}{4\pi \,\mathcal{L}_i}}
\ e^{-3\kappa_i^2/4\,\mathcal{L}_i}
+O\left(\frac{1}{\varepsilon\sqrt{n}}\right)
\right),\\
&=&
\ \parth{\frac 9{8n}}^{\delta}\xi(\ell,k)
\ \ \left(1+\mathcal O_K\left(n^{-1/2}\right)\right),
\end{eqnarray*}
in which the sequence of normalised variations of
      labels along superedges, $k=(k_i)_{1\leq i\leq2p-1}$
is the inverse image of
$(w,\mu)$ through $M_T$.
In other terms, $$(\ell,k)=\tilde M_T^{-1}(x,m,w,\mu),$$
and the previous relation can be written
\begin{eqnarray*}
\lefteqn{\mathbb{P}\left(Z^{(n)}(\tau)=(x,m,w,\mu)\mid
Y^{(n)}=(x,m)\right)}\\
&=&\parth{8n/9}^{-\delta}
\, \xi_{T_\sigma}(x,m,w,\mu)
\, \big(1+\mathcal O_K\big(n^{-1/2}\big)\big),
\end{eqnarray*}
leading to the desired result, through  Proposition~\ref{pro:fdd-shape}.
\end{proof}

\begin{proof}[Proof of Corollary~\ref{cor:weakZ}]
In view of the measure preserving property of $\tilde M_T$, we have
\[
\mathbb{E}(\Phi(Z(\tau)))\,=
\int\Phi\circ \tilde M_T(\ell,k)\;f\circ \tilde M_T(\ell,k)\,d\ell\,dk\,.
\]
Therefore,  the lemma follows upon proving that
\[
\lim_n\left|\mathbb{E}\left(\Phi(Z^{(n)}(\tau))\right)-
\int\Phi\circ \tilde M_T(\ell,k)\;f\circ \tilde M_T(\ell,k)\,d\ell\,dk\,\right|
=0.\]
Set
\begin{equation*}
\phi_T=\Phi\circ \tilde M_T,\hspace{15mm}f_T=f\circ \tilde M_T.
\end{equation*}
The function $\phi_T$ has a compact support $\tilde
K=\tilde M_T^{-1}K$ that is included in $\tilde
M_T^{-1}K_{\sigma,\epsilon}
\subset(\epsilon,\infty)^{2p-1}\times\mathbb{R}^{2p-1}$.  Since
$\tilde K$ is compact, there exists an $\epsilon'>0$ and a compact
$K'$ such that $\tilde K\subset
K'\subset(\epsilon,\infty)^{2p-1}\times\mathbb{R}^{2p-1}$ and
$d(\partial \tilde K,\partial K')>\epsilon'$. We shall use the fat
boundary $K'\setminus\tilde K$, on which $\phi_T$ is identically zero,
to deal with boundary effects.

Let finally $\tilde\Delta(n,K)=\tilde M_T^{-1}\Delta(n,K)$ be the
discretized version of $\tilde K$ and similarly $\Delta'(n,K)=\tilde
M_T^{-1}\Delta(n,M_TK')$ that of $K'$; by construction
$\tilde\Delta(n,K)\subset\Delta'(n,K)$ are finite  sets with
$\mathcal O_K(n^{3\delta})$ elements,
and $\phi_T$ is identically
zero over $\Delta'(n,K)\setminus\tilde\Delta(n,K)$.

First, we have
\begin{eqnarray*}
{\mathbb{E}\left(\Phi(Z^{(n)}(\tau))\right)}
&=&
\sum_{(\ell,k)\in\Delta'(n,K)}\phi_T(\ell,k)
\,\mathbb{P}\left((\ell^{(n)},k^{(n)})=(\ell,k)\right)\\
&=&
(2n)^{-2\delta}(8n/9)^{-\delta}
\!\!\sum_{(\ell,k)\in\Delta'(n,K)}\!\!\phi_T(\ell,k)\,f_T(\ell,k)\\
&&\;\;+\;\;\Vert\Phi\Vert_\infty\cdot\mathcal O_K(n^{-1/2}),
\end{eqnarray*}
the second equality due to the local limit convergence
(Lemma~\ref{lem:locZ}).  Next, difference between the discrete
summation and the integral is bounded in terms of the modulus of
continuity $\omega(\phi_T\cdot f_T,K',\cdot)$ of $\phi_T\cdot f_T$ on
$\tilde K$ (recall that the modulus of continuity of a function $g$
on a compact $K$ is
$\omega(g,K,\epsilon)=\sup_{0<d(x,y)<\epsilon}|g(x)-g(y)|$, and that
if $g$ is uniformly continuous on $K$, it satisfies
$\omega(g,K,\epsilon)=O(\epsilon)$ as $\epsilon$ tends to zero). This
yields
\begin{eqnarray*}
&&\left|
(2n)^{-2\delta}(8n/9)^{-\delta}
\!\!\sum_{(\ell,k)\in\Delta'(n,K)}\!\!\phi_T(\ell,k)\,f_T(\ell,k)
-
\int\phi_T(\ell,k)\,f_T(\ell,k)\,d\ell\,dk
\right|\\
&&\qquad\leq\;\mathrm{measure}(K')\;\cdot\;
\omega(\phi_T\cdot f_T,K', n^{-1/4}).
\end{eqnarray*}
Observe that in this summation the compact $K'$ has been approximated
by a union of boxes of diameter $O(n^{-1/4})$ and boundary effect
should be considered. However $\phi_T$ is identically
zero on a region $K'\setminus K$ with $d(\delta K,\delta
K')>\varepsilon'$, that contains all boxes intersecting
the boundary for $n$ large enough. The boundary effect is thus null.
\end{proof}

\subsubsection{Weak convergence of $Z^{(n)}(\tau)$.}

According to the Porte-Manteau
Theorem \cite[Ch. 1]{BILL}, we need
$$\lim_n \mathbb{E}\ \left[\Phi\parth{Z^{(n)}(\tau)}\right]=\int\Phi\
f =\mathbb{E}\croc{\Phi\parth{Z(\tau)}}$$ to hold for any bounded
uniformly continuous $\Phi$.  Now consider
\[K_m\ =\ \overline B(0,\rho_m)\bigcap\parth{\build{\bigcup}{\sigma}{}
R_{\sigma,\varepsilon_m}\times\mathbb{R}^{2p-1}},\] in which
$\overline B(0,\rho_m)$ is the closed ball, in $\mathbb R^{4p-2}$,
with radius $\rho_m$ and let simultaneously $\rho_m$ increase to
$+\infty$ and $\varepsilon_m$ decrease to $0$.  We obtain a increasing
sequence of compacts $K_m$, each of these compacts having $(p-1)!$
connected components.  As the limit of this sequence has a
Lebesgue--negligible complement in $\mathbb R^{4p-2}$, we can choose
the sequences $(\rho_m)_{m>1}$ and $(\varepsilon_m)_{m>1}$ in such a
way that
$$\mathbb{P}\parth{Z(\tau)\in K_m}\ge 1-\frac 1m.$$
There exist uniformly continuous functions $\Psi_m:\ \ \mathbb{R}^{4p-2}
\longrightarrow [0,1]$ ,  such that
\[\Psi_m\vert_{K_m}
\equiv 1, \hspace{1cm}\Psi_m\vert_{K_{m+1}^c}
\equiv 0. \]
By construction,
\[
\left|
\mathbb{E}\Big(\Phi(Z(\tau))\Big)\;-\;
\mathbb{E}\Big(\Psi_m\cdot\Phi\;(Z(\tau))\Big)
\right|
\;\leq\; \frac{\Vert\Phi\Vert_\infty}m.,
\]
Moreover the product $\Psi_m\cdot\Phi$ is now a finite sum of functions
satisfying the assumptions of Corollary~\ref{cor:weakZ}: this yields
\begin{eqnarray*}
\lim_n\  \mathbb{E}\left(\Psi_m\cdot\Phi\;(Z^{(n)}(\tau))\right)&=&\mathbb{E}\
\left(\Psi_m\cdot\Phi\;(Z(\tau))\right).
\end{eqnarray*}
Next observe that, by definition of $\Psi_m$,
\[
\mathbb{P}\left(Z^{(n)}(\tau)\in K_m\right)\;\geq\;
\mathbb{E}\left(\Psi_{m-1}(Z^{(n)}(\tau))\right),
\]
and
\[
\mathbb{E}\big(\Psi_{m-1}(Z(\tau))\big)\;\geq\;
\mathbb{P}\big(Z(\tau)\in K_{m-1}\big).
\]
Therefore, applying Corollary~\ref{cor:weakZ} to $\Psi_{m-1}$,
\begin{eqnarray*}
\liminf_n\  \mathbb{P}\left(Z^{(n)}(\tau)\in K_m\right)&\geq&
\mathbb{P}\left(Z(\tau)\in K_{m-1}\right)\;\ge\; 1-\frac 1{m-1}.
\end{eqnarray*}
Moreover
\begin{eqnarray*}
\mathbb{E}\left((1-\Psi_m)\cdot\Phi\;(Z^{(n)}(\tau))\right)
&\leq& \mathbb{P}\left(\Psi_m(Z^{(n)}(\tau))<1\right)\cdot
\Vert\Phi\Vert_\infty,\\
&\leq& \mathbb{P}\left(Z^{(n)}(\tau)\not\in K_m\right)\cdot
\Vert\Phi\Vert_\infty.
\end{eqnarray*}
Thus, taking limit for $n$ going to infinity and applying the previous
lower bound,
\begin{eqnarray*}
\limsup_n\mathbb{E}\left((1-\Psi_m)\cdot\Phi\;(Z^{(n)}(\tau))\right)
&\leq&\frac{\Vert\Phi\Vert_\infty}{m-1}.
\end{eqnarray*}
Finally, the decomposition
\begin{eqnarray*}
\mathbb{E}\left(\Phi(Z^{(n)}(\tau))\right)
&=&
\mathbb{E}\left(\Psi_m\cdot\Phi\;(Z^{(n)}(\tau))\right)
\;+\;
\mathbb{E}\left((1-\Psi_m)\cdot\Phi\;(Z^{(n)}(\tau))\right),
\end{eqnarray*}
yields
\begin{eqnarray*}
\limsup_n\left|\mathbb{E}\left(\Phi(Z^{(n)}(\tau))\right)
-\mathbb{E}\big(\Phi(Z(\tau))\big)\right|
\leq \frac{2\Vert\Phi\Vert_\infty}{m-1}.
\end{eqnarray*}
Letting $m$ go to infinity gives the weak convergence of
$Z^{(n)}(\tau)$ to $Z(\tau)$ as claimed. The convergence of
$X^{(n)}(\tau)$ is a by-product.

\subsection{A deviation bound for the largest label}
\label{sec:bound}

In this section, a rough but exponential deviation bound for the value
of the largest label in a forest of $k$ embedded trees with $n$ edges
is obtained. For $k=1$, Proposition~\ref{pro:devbound} is exactly obtained.

Let $\mathcal{EV}_{k,2n}$ denote the set of $k$-uples of element of
$\mathcal{EV}$ of total length $2n$: an element
$[(E_1,V_1),\ldots,(E_k,V_k)]$ of $\mathcal{EV}_{k,2n}$ codes for a
forest of $k$ embedded trees (each $(E_i,V_i)$ codes for a tree,
according to Proposition~\ref{pro:contours}). Equivalently, one may
concatenate the $k$ pairs and view any element of
$\mathcal{EV}_{k,2n}$ as a pair $(E,V)=(E_1\cdots E_k,V_1\cdots
V_k)\in\mathcal{EV}_{2n}$ together with a set of concatenation times
$0=t_0\leq t_1\leq\cdots\leq t_k=2n$, subject to the conditions
$E(t_i)=V(t_i)=0$ for all $i=1,\ldots,k$. In this identification,
$E_i$ has length $2n_i=t_i-t_{i-1}$.

Let $(E^{k,n},V^{k,n})$ denote a random forest of
$\mathcal{EV}_{k,2n}$ under the uniform distribution.
\begin{pro}\label{pro:loosebound}
There exists $y_0>0$ such that, for all $n$, $k$ and $y\ge y_0$,
\begin{equation}\label{equ:loosebound}
{\textstyle\mathbb{P}}
\left(\sup_{0\leq t\leq 2n}V^{k,n}(t)
\,>\, yn^{1/4}\right)\;<\;e^{-y}.
\end{equation}
\end{pro}
This is the key fact in the proof of tightness of the sequence
        $X^{(n)}$, given in the last subsection, and it also leads to the
        convergence of moments in Corollary \ref{radcv} through the
        following weaker formulation (with $k=1$): for $y>y_0$, and for
        all~$n$,
\begin{equation}
\mathbb{P}\left(\sup_{0\leq s\leq1}\hat
W^{(n)}_s\,>\,(8/9)^{1/4}y\right)\;\leq e^{-y},
\end{equation}
which is exactly Proposition~\ref{pro:devbound}.

The proof is based on a \emph{branch decomposition}, that is discussed
in the next paragraph. Then, after two preliminary results on
parameters of the middle branch of a random tree
(Paragraphs~\ref{sec:middlebranchlength}
and~\ref{sec:middlebranchlargest}), Proposition~\ref{pro:loosebound} is
proved by induction (Paragraphs~\ref{sec:induct3}, \ref{sec:induct4}
and~\ref{sec:induct5}).  At the price of more technical details in
these latter paragraphs, the bound could be improved to
$e^{-c_\epsilon y^{4/3-\epsilon}}$ for any fixed $\epsilon>0$.

\subsubsection{The branch decomposition at time $t$}\label{subs:brandw}
Let $(E,V)\in\mathcal{EV}_{k,2n}$ with concatenation times $0\leq
t_1\leq\ldots\leq t_k\leq2n$ and let $t\in(0,2n)$.  Suppose moreover
that $t_{p-1}<t<t_p$, that is $t$ occurs during the contour traversal of
$(E_p,V_p)$, the $p$th component of the forest encoded by $(E,V)$ and
let $U$ be the tree encoded by $(E_p,V_p)$.

To any vertex $v$ of $U$ is associated the set of edges in the unique
simple path from $v$ to the root of $U$, denoted $\mathrm{br}(v)$ (for
the \emph{branch} of $v$).
\begin{figure}
\begin{center}
\includegraphics[width=11cm]{\rep 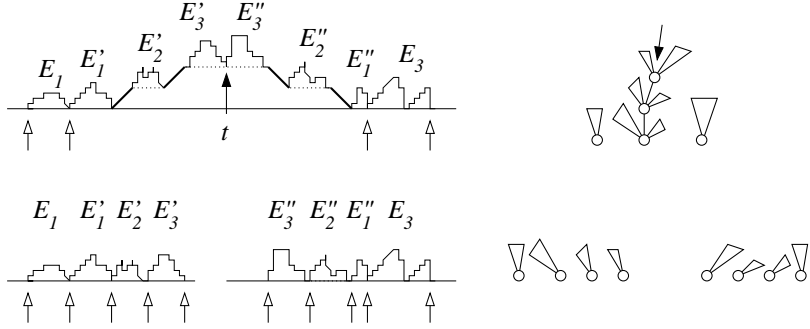}
\end{center}
\caption{Decomposition of a forest at time $t$. Arrows under the walks
indicate concatenation times.
\label{fig:branch}}
\end{figure}
Recall that $v_t$ denote the vertex visited
at time $t$ of the contour traversal of $U$.  Observe that the height
$E(t)$ is $|\textrm{br}(v_t)|$, the length of the branch from
the root $r$
to $v_t$, while the label $V(t)$ of $v_t$ is given by
\begin{equation}\label{iidonbranches}
V(t)=\sum_{\epsilon\in\textrm{br}(v_t)}\kappa(\epsilon).
\end{equation}
Let
$\ell=E(t)=|\mathrm{br}(v_t)|$, and call $\epsilon_i$ the edge of
$\mathrm{br}(v_t)$ between heights $i-1$ and $i$, for $i=1,\ldots,\ell$.
The maximal label on the branch $\textrm{br}(v_t)$ is
\[
\bar H(t)\;=\;\sup_{0\leq j\leq \ell}\;\sum_{i=1}^j \kappa(\varepsilon_i)
\]
where $\varepsilon_i$ is the $i$th edge of $\textrm{br}(v_t)$.

The branch of $v_t$ induces a decomposition of $U$ into two forests of
trees,   the \emph{branch decomposition}, which we now
phrase in terms of $(E,V)$.
\begin{itemize}
\item From time $t_{p-1}$ to $t$, the edges $\epsilon_1$, \ldots,
$\epsilon_{\ell}$ are successively traversed away from the root on the
branch $\mathrm{br}(v_t)$, at times
$t'_1<\cdots<t'_{\ell}<t$.
Let $(E'_i,V'_i)$ be the part of the contour walks $(E,V)$ between
        $t'_i$ and $t'_{i+1}$ (with the convention
that $t'_0=t_{p-1}$ and $t'_{\ell+1}=t$).
\item From time $t$ to $t_p$, the edges $\epsilon_\ell$, \ldots,
$\epsilon_1$  are successively traversed back toward the root on the
branch $\mathrm{br}(v_t)$, at times
$t''_{\ell+1}<\cdots<t''_{1}$.
Let $(E''_i,V''_i)$ be the part of the contour walks $(E,V)$ between
        $t''_{i+1}$ and $t''_i$ (with the convention that
$t''_{\ell+1}=t$ and $t''_0=t_{p}$).
\end{itemize}
The contour walks $(E'_i,V'_i)$ (resp. $(E''_i,V''_i)$) encodes for the
left (resp. right) subtree attached at the $i$th vertex of the branch
$\textrm{br}(v_t)$.
One can see $t'_i$ as the time of the last upcrossing of heights $(i-1,i)$
         before time $t$, and $t''_i$ as the time of the first downcrossing of
         heights $(i,i-1)$ after time $t$.

Upon shifting the walks $E'_i$ and $E''_i$ down by $i$ so that they
start from zero, and also shifting the walks $V'_i$ and $V''_i$ by
$V(t'_i)$ (resp. $V(t''_{i+1})$): the two forests
\begin{eqnarray*}
\textrm{left}(t)&=&[(E_1,V_1),\ldots,(E_{p-1},V_{p-1}),
(E'_0,V'_0),\ldots,(E'_\ell,V'_\ell)]
\\
\textrm{right}(t)&=&[(E''_\ell,V''_\ell),\ldots,(E''_0,V''_0),
(E_{p+1},V_{p+1}),\ldots,(E_k,V_k)]
\end{eqnarray*}
belong respectively to $\mathcal{EV}_{k',n'}$ and
$\mathcal{EV}_{k'',n''}$ with $k'=p+\ell$, $n'=t-\ell$,
$k''=k-p+1+\ell$, and $n''=2n-t-\ell$.

Let us apply this branch decomposition to a random forest
$(E^{k,n},V^{k,n})$.
In view of expression~(\ref{iidonbranches}), conditionally on
$E^{(k,n)}(t)=\ell$,
\[
V^{(k,n)}(t)\;\mathop{=}^{\textrm{law}}\;S_\ell,
\]
where $S_k$ denote the sum of $k$ i.i.d. random variables uniform on
$\{-1,0,1\}$. Moreover, again conditionally on $E^{(k,n)}(t)=\ell$,
\[
\textrm{left}^{(k,n)}(t)
\;\mathop{=}^{\textrm{law}}\;(E^{k',n'},V^{k',n'}),
\qquad
\textrm{right}^{(k,n)}(t)
\;\mathop{=}^{\textrm{law}}\;(E^{k'',n''},V^{k'',n''}).
\]

\subsubsection{The middle branch length in a random forest}
\label{sec:middlebranchlength}
      The first step
of the proof is a bound for the tail probability of $E^{k,n}(n)$, the
length of the middle branch ($t=n$).
\begin{lem}  For $4\ell^2>27\cdot n$,
\begin{eqnarray*}
\Pr(E^{k,n}(n)>\ell)&<&A\cdot\frac{\ell}{n^{1/2}}
\exp\left({-\frac{4\ell^2}{27n}}\right).
\end{eqnarray*}
\end{lem}
\begin{proof}
With the notations of Subsection \ref{subs:brandw}, for $1\le p\le k$,
let $t^{k,n}_p$ denote the $p$th concatenation time of $E^{k,n}$. Then
\[\Pr\parth{E^{k,n}(n)=\ell\,\left\vert\,
t^{k,n}_{p+1}-t^{k,n}_p=2m,\, 0\le n-t^{k,n}_p=a\le 2m \right.}
=\Pr\parth{E^{(m)}(a)=\ell},\]
so the proof reduces to bound  $\Pr\parth{E^{(m)}(a)>\ell}$  uniformly
on pairs $(a,m)$ such that $0\le a\le 2m\le 2n$. We have
\begin{eqnarray*}
\Pr\parth{E^{(m)}(a)=\ell}&=&\frac{C(a;\ell)\,C(2m-a;\ell)}{C(2m;0)}\\
&=&\frac{C(\alpha;\ell)\,C(\beta;\ell)}{C(2m;0)},
\end{eqnarray*}
with an obvious change of variables.
       From \cite[Ch. I.3, Th~2]{bolobas}, for $k>n/2+h$ and $1\leq h\leq
n/6$
\[
{\sqrt{\frac{\pi n}2}}\binom{n}{k}\;<\;2^n
\exp\left({-\frac{2h^2}{n}+\frac{2h}n+\frac{4h^3}{n^2}}\right).
\]
With $n=\alpha+1$, $k=(\alpha+\ell)/2+1$, $1\leq
h=\ell/6\leq \alpha/6\leq n/6$,  this bound yields:
\begin{eqnarray*}
C(\alpha;\ell)&=&\frac{\ell+1}{\alpha+1}{\alpha+1\choose (\alpha+\ell)/2+1}\\
&<&\sqrt{\frac 8{\pi}}\cdot\frac\ell{n^{3/2}}\;2^n
\exp\left({-\frac{\ell^2}{18n}
+\frac{\ell}{3n}+\frac{\ell^3}{54n^2}}\right),\\
&<&e^{1/3}\sqrt{\frac 8{\pi}}\cdot\frac\ell{n^{3/2}}\;2^n
\exp\left({-\frac{\ell^2}{27n}}\right)\\
&<&2e^{10/27}\sqrt{\frac 8{\pi}}\cdot\frac\ell{\alpha^{3/2}}\;2^{\alpha}
\exp\left({-\frac{\ell^2}{27\alpha}}\right),
\end{eqnarray*}
Thus,
\begin{eqnarray*}
\Pr\parth{E^{(m)}(a)=\ell}
&<&B_1\cdot\frac{\ell^{2}m^{3/2}}{(\alpha\beta)^{3/2}}
\exp\left({-\frac{\ell^2}{27\alpha}
{-\frac{\ell^2}{27\beta}}}\right),\\
&<&B_2\cdot\frac{\ell^{2}}{m^{3/2}}
\exp\left({-\frac{4\ell^2}{27m}}\right),
\end{eqnarray*}
the latter inequality since the maximum
of the function $(xy)^{-3/2}\exp\parth{-x^{-1}-y^{-1}}$, subject to
$x+y=27m/\ell^2<1$, is obtained for $x=y$.  The last inequality entails that,
for $4\ell^2>27m$,
\begin{eqnarray*}
\Pr\parth{E^{(m)}(a)>\ell}
&<&B_3\cdot\frac{\ell}{m^{1/2}}
\exp\left({-\frac{4\ell^2}{27m}}\right)\\
&<&B_3\cdot\frac{\ell}{n^{1/2}}
\exp\left({-\frac{4\ell^2}{27n}}\right).
\end{eqnarray*}
\end{proof}

\subsubsection{The largest label on the middle branch}
\label{sec:middlebranchlargest}
\begin{lem}\label{lem:branch}
Let $\bar H^{k,n}$ be the largest label on the branch between the root
and the vertex reached at $t=n$ by the contour traversal of
$(E^{k,n},V^{k,n})$. Then there exists $c_0$  such that
for all $k$, $n$ and $h$,
\begin{equation}
\mathbb{P}(\bar H^{k,n}>h)\leq c_0\ e^{-\frac1{9}(hn^{-1/4})^{4/3}}.
\end{equation}
\end{lem}
\begin{proof}
As already discussed the conditional probability that the largest
label is $h$ knowing that the branch has length $\ell$ is exactly the
probability that a random walk with steps $\{1,0,-1\}$ of length
$\ell$ has maximal value $h$. Using the reflection principle, Azuma's
inequality
         \cite[Th.~2.1, p. 85]{AS} then reads
\begin{eqnarray*}
\mathbb{P}(\bar H^{k,n}>h\mid E^{k,n}(n)=\ell)&\leq&
         2\ e^{-h^2/2\ell}.
\end{eqnarray*}
Next, as previously calculated, for $4\ell^2\geq 27n$,
\begin{eqnarray*}
\mathbb{P}(E^{k,n}(n)>\ell)\;<\;A \cdot(\ell/n^{1/2})\;e^{-4\ell^2/27n}
\end{eqnarray*}
Finally, assuming $h>(27/4)^{3/4}n^{1/4}$, so that the previous
inequality holds,
\begin{eqnarray*}
\mathbb{P}(\bar H^{k,n}>h)
&\leq&\mathbb{P}(\bar H^{k,n}>h\mid E^{k,n}(n)\leq h^{2/3}n^{1/3})\;+\;
\mathbb{P}(E^{k,n}(n)>h^{2/3}n^{1/3})\\
&\leq&
2 e^{-h^{4/3}/2n^{1/3}}\;+\;A\cdot(h^{2/3}n^{-1/6})\,e^{-4h^{4/3}/27n^{1/3}}.
\end{eqnarray*}
The lemma follows for  $h>(27/4)^{3/4}n^{1/4}$, taking
       $c_0$ large enough, and also  for
       $h\le(27/4)^{3/4}n^{1/4}$, taking
       $c_0\ge e^{3/4}$.
\end{proof}

\subsubsection{Conditional induction in the case $E(n)>0$}
\label{sec:induct3}
\begin{lem}\label{lem:induct1}
Assume the bound (\ref{equ:loosebound}) holds true for some $y_0$
for $V^{k,m}$ with $m<n$. Then, for all $\ell>0$ and $k\geq p\geq0$,
the probability
\begin{eqnarray*}
\lefteqn{p_{k,n}(y;h,\ell,p)\;=}\\
&&\mathbb{P}\Big(\sup_{0\leq t\leq n}V^{k,n}(t)\;>\;yn^{1/4}\;\Big|\;
\bar H^{k,n}=h,\, E^{k,n}(n)=\ell,\, t^{k,n}_{p-1}<n<t^{k,n}_{p}\,\Big),
\end{eqnarray*}
satisfies
\begin{equation}\label{ind-bound}
p_{k,n}(y;h,\ell,p)\leq 4e^{-2^{1/4}(y-hn^{-1/4})},
\quad\textrm{ provided $h\le n^{1/4}(y-y_0/2^{1/4})$.}
\end{equation}
\end{lem}
\begin{proof}
Assume $(E^{k,n},V^{k,n})=(E_1\cdots E_k,V_1\cdots V_k)$ is such that
$E^{k,n}(n)=\ell>0$ and $t^{k,n}_{p-1}<n<t^{k,n}_{p}$, so that
$t=n$ occurs inside $(E_p,V_p)$. Apply the branch
decomposition (see Section~\ref{subs:brandw}) at $t=n$ to
$(E^{k,n},V^{k,n})$ and let
\begin{eqnarray*}
(\bar E^{k,n},\bar V^{k,n})&=&\textrm{left}^{k,n}(n)
\;=\;[(E_1,V_1),\ldots,(E_{p-1},V_{p-1}),
(E'_0,V'_0),\ldots,(E'_\ell,V'_\ell)],\\
(\bar{\bar E}^{k,n},\bar{\bar V}^{k,n})&=&\textrm{right}^{k,n}(n)
\;=\;[(E''_\ell,V''_\ell),\ldots,(E''_0,V''_0),
(E_{p+1},V_{p+1}),\ldots,(E_k,V_k)].
\end{eqnarray*}
Upon taking $k'=p+\ell$, $k''=k-p+\ell+1$ and $n'=(n-\ell)/2$,
\[
(\bar E^{k,n},\bar V^{k,n})\;\mathop{=}^{\textrm{law}}\;
(E^{k',n'},V^{k',n'})
\quad\textrm{ and }\quad
(\bar{\bar E}^{k,n},\bar{\bar V}^{k,n})\;\mathop{=}^{\textrm{law}}\;
(E^{k'',n'},V^{k'',n'})
\]
Observe now that, in the previous decomposition,
\[
\sup_{0\leq t\leq 2n}V(t)\;\leq\;
\sup\Big(\bar H+\sup_{0\leq t\leq 2n'}\bar V(t)\,,\;
\bar H+\sup_{0\leq t\leq 2n'}\bar{\bar V}(t)\Big),
\]
so that
\begin{eqnarray*}
\lefteqn{p_{k,n}(y;h,\ell,p)\leq}\\
&&\mathbb{P}\Big(
\sup_{0\leq t\leq 2n'}\bar V^{k,n}(t)>yn^{1/4}-h\;
\Big|\;\bar H^{k,n}=h,\, E^{k,n}(n)=\ell,\,
t^{k,n}_{p-1}<n<t^{k,n}_{p}\,\Big)\\
&+&
\mathbb{P}\Big(
\sup_{0\leq t\leq 2n'}\bar{\bar V}^{k,n}(t)>yn^{1/4}-h\;
\Big|\;\bar H^{k,n}=h,\, E^{k,n}(n)=\ell,\,
t^{k,n}_{p-1}<n<t^{k,n}_{p}\,\Big).
\end{eqnarray*}
Hence, in view of the preceding identities in law,
\begin{eqnarray*}
p_{k,n}(y;h,\ell,p)&\leq&
\mathbb{P}\Big(
\sup_{0\leq t\leq 2n'}V^{k',n'}(t)>yn^{1/4}-h\;\Big)\\
&&\!\!\!\!\!+\;\;
\mathbb{P}\Big(
\sup_{0\leq t\leq 2n'}V^{k'',n'}(t)>yn^{1/4}-h\;\Big).
\end{eqnarray*}
Observe that $2n'\le n$. Hence
\begin{eqnarray*}
p_{k,n}(y;h,\ell,p)&\leq&
\mathbb{P}\Big(
\sup_{0\leq t\leq 2n'}V^{k',n'}(t)>2^{1/4}(y-hn^{-1/4})n'{}^{1/4}\;\Big)\\
&&\!\!\!\!\!+\;\;
\mathbb{P}\Big(
\sup_{0\leq t\leq 2n'}V^{k'',n'}(t)>2^{1/4}(y-hn^{-1/4})n'{}^{1/4}\;\Big).
\end{eqnarray*}
The induction hypothesis now implies, for $2^{1/4}(y-hn^{-1/4})\ge y_0$,
that is, for all $h\le n^{1/4}(y-y_0/2^{1/4})$,
\begin{equation*}
p_{k,n}(y;h,\ell,p)\leq 2e^{-2^{1/4}(y-hn^{-1/4})},
\end{equation*}
which is exactly the lemma, up to a factor 2 added for later
convenience.
\end{proof}
\subsubsection{Conditional induction in the case $E(n)=0$}\label{sec:induct4}
\begin{lem}\label{lem:induct2}
Assume the bound (\ref{equ:loosebound}) holds true for some $y_0$ for
$V^{k,m}$ with $m<n$. Then, provided $y\geq y_0/2^{1/4}$,
\begin{eqnarray*}
\mathbb{P}\Big(\sup_{0\leq t\leq 2n}V^{k,n}(t)\;>\;yn^{1/4}\;\Big|\;
E^{k,n}(n)=0\,\Big)\;\leq\; 4e^{-2^{1/4}y}.
\end{eqnarray*}
The bound~(\ref{ind-bound}) thus remain valid in the case $\ell=0$.
\end{lem}
\begin{proof}
In this case the decomposition at $t=n$ is even simpler.  Let
$p=\sup\{i\mid t_i< n\}$, and $q=\inf\{i\mid t_{i}> n\}$,
and consider a decomposition in four parts, cutting at times $t_p$,
$n$ and $t_q$.  The two contour walks for $0\leq t\leq t_p$ and
$t_q\leq t\leq 2n$ are uniform on $\mathcal{EV}_{p,t_p}$ and
$\mathcal{EV}_{k-q,2n-t_q}$.  The other two contour walks are uniform
on $\mathcal{EV}_{n-t_p}$ and $\mathcal{EV}_{t_q-n}$.

The result then follows using the induction hypothesis on the 4  parts.
\end{proof}\smallskip

\subsubsection{Complete induction and proof of
Proposition~\ref{pro:loosebound}}\label{sec:induct5}~\smallskip

Observe that the bounds of Lemmas~\ref{lem:induct1} and
\ref{lem:induct2} do not depend on $\ell$ or $p$, so that,
assuming (\ref{equ:loosebound}) holds for some $y_0$
for $V^{k,m}$ with $m<n$,
\begin{eqnarray*}
\mathbb{P}\Big(\sup_{0\leq t\leq 2n}V^{k,n}(t)\;>\;yn^{1/4}\;\Big|\;
\bar H^{k,n}=h\,\Big)
&\leq& 4e^{-2^{1/4}(y-hn^{-1/4})},
\end{eqnarray*}
provided $h\le n^{1/4}(y-y_0/2^{1/4})=h_0$.
Using this bound,
\begin{eqnarray*}
f_{k,n}(y)&=&
{\mathbb{P}}
\Big(\sup_{0\leq t\leq 2n}V^{k,n}(t)\,>\,yn^{1/4}\Big)
\\&\leq&
4e^{-2^{1/4}y}\sum_{h=0}^{h_0}
e^{2^{1/4}hn^{-1/4}}\mathbb{P}(\bar H^{k,n}=h)
\;+\;\mathbb{P}\parth{\bar H^{k,n}>h_0},\\
&\leq&4e^{-2^{1/4}y}\parth{1+\sum_{h=1}^\infty
(e^{2^{1/4}hn^{-1/4}}-e^{2^{1/4}(h-1)n^{-1/4}})
\mathbb{P}(\bar H^{k,n}\geq h)}\\
&&+\;\mathbb{P}(\bar H^{k,n}>{n^{1/4}(y-y_0/2^{1/4})}).
\end{eqnarray*}
In view of Lemma~\ref{lem:branch}, the summation is bounded
by a convergent integral  that evaluates to a
constant $c_1$.
Lemma~\ref{lem:branch} allows also to dispose of the second
term:
\begin{eqnarray*}
f_{k,n}(y)
&\leq& 4e^{-2^{1/4}y}(1+c_1)+c_0e^{-\frac1{9}(y-y_0/2^{1/4})^{4/3}}
\end{eqnarray*}
Now observe that
\begin{eqnarray*}
e^{-\frac1{9}(y-y_0/2^{1/4})^{4/3}}
&=&
e^{-y}\ e^{y_0/2^{1/4}}\ e^{(y-y_0/2^{1/4})
-\frac1{9}(y-y_0/2^{1/4})^{4/3}}\\
&\leq&
e^{-y}\ e^{y_0/2^{1/4}}\ e^{(y_0-y_0/2^{1/4})
-\frac1{9}(y_0-y_0/2^{1/4})^{4/3}}\\
&=&
e^{-y}\ e^{y_0 -\frac1{9}(1-1/2^{1/4})^{4/3}y_0^{4/3}}
\end{eqnarray*}
for $y\ge y_0$, as soon as $x\rightarrow x
-\frac1{9}x^{4/3}$ is decreasing on the interval
       $\left[y_0-y_0/2^{1/4}, +\infty\right[$, that is, for $y_0\ge 1933$.
The bound can thus be rewritten as
\begin{eqnarray*}
f_{k,n}(y)
&\leq&e^{-y}\left(4e^{-(2^{1/4}-1)y_0}(1+c_1)+
c_0e^{y_0 -\frac19(1-1/2^{1/4})^{4/3}y_0^{4/3}}
\right)
\end{eqnarray*}
which is smaller than $e^{-y}$ for $y_0$ large enough, so that
induction can be carried on (Recall that $c_0$ and $c_1$ do not
depend on $n$).
The case $n=1$ holds true for $y_0\ge 1$, and the proof of
Proposition~\ref{pro:loosebound} is complete.

\subsection{Tightness}\label{sec:tight}
Tightness for  the bidimensional path follows from the tightness
          of the two projections. The tightness of the first
          projection was proved by
Kaigh \cite{KAIGH}.
Thus we  only have to prove the following proposition (cf. \cite[Ch. 3]{BILL}).
\begin{pro}
For all $\varepsilon>0$ and $\delta>0$ there exists $m$ such that for $n$
          large enough
\[\mathbb{P}\parth{\sup_{0\le s\le 1}\left\vert \hat W^{(n)}_s
-\hat W^{(n)}_{\lfloor ms\rfloor/m} \right\vert\,\ge\,\delta}\le
     \varepsilon.\]
\end{pro}
We proceed by bounding, for $m$ large enough and all $i$ with $1\leq
i\leq m$,
\begin{eqnarray*}
p_{m,i}(n) &=& \mathbb{P}\parth{\sup_{(i-1)/m\le s\le i/m}\left\vert
\hat W^{(n)}_s
-\hat W^{(n)}_{(i-1)/m} \right\vert\,\ge\,\delta}\\
&=&
\mathbb{P}\parth{\sup_{t_1\leq t\leq t_2}\left\vert
V^{(n)}(t)-V^{(n)}(t_1)\right\vert\,\ge\,\delta'n^{1/4}},
\end{eqnarray*}
where $t_1=\lfloor 2(i-1)n/m\rfloor$, $t_2=\lfloor 2in/m\rfloor$ and
$\delta'=(8/9)^{1/4}\delta$.

The proposition is an immediate consequence of the following lemma,
upon taking $m$ large enough and summing on $i$.
\begin{lem}
For all $\varepsilon>0$ and $\delta>0$ there exist $m_0$ such that for
all $m>m_0$ and $1\leq i\leq m$,
\[
\exists n_0,\;\forall n>n_0,\qquad p_{m,i}(n)\;\leq\;
\frac{\varepsilon}{2m}+C_\delta(m),
\]
where $C_\delta(m)$ is exponentially decreasing
in a positive power of $m$.
\end{lem}
\begin{proof}
For simplicity of notations let us assume $t_2-t_1=2n/m$. (The general
case is identical but obscured by a collection of $\lfloor n/m\rfloor$.)

Let us consider the shape $T$ associated with the times $t_1$ and
$t_2$, as defined in Section~\ref{subs:shapes}. The branches
$\textrm{br}(v_{t_1})$ from $v_{t_1}$ to the root and
$\textrm{br}(v_{t_2})$ from $v_{t_2}$ to the root meet at the unique
branchpoint $v$ of the shape $T$. Let
$B=\textrm{br}(v_{t_1})\setminus\textrm{br}(v_{t_2})$ be the branch
between $v_{t_1}$ and $v$ and
$C=\textrm{br}(v_{t_2})\setminus\textrm{br}(v_{t_1})$ the branch
between $v_{t_2}$ and $v$. Between $t_1$ and $t_2$, a total of $k$
edges $\epsilon_i$, $i=1,\ldots,k$ of $U$ are successively traversed
on the branch from $v_{t_1}$ to $v_{t_2}$ by the contour walk. The
$k'$ first are along the branch $B$ and traversed toward the root,
while the $k''$ next are along the branch $C$ and traversed away from
the root. By construction,
\begin{eqnarray*}
k'&=&E^{(n)}(t_1)-\inf_{[t_1,t_2]}E^{(n)}
\;=\;\sqrt{2n}\ \croc{e^{(n)}
\parth{\frac {i-1}{m}}-\inf_{[\frac{i-1}{m},\frac{i}m]}e^{(n)}
},\\
k''&=&E^{(n)}(t_2)-\inf_{[t_1,t_2]}E^{(n)}
\;=\;\sqrt{2n}\ \croc{e^{(n)}
\parth{\frac {i}{m}}-\inf_{[\frac{i-1}{m},\frac{i}m]}e^{(n)}
},
\end{eqnarray*}
and the total length $\Delta^{(n)}_i=k=k'+k''$ of the branch from
$v_{t_1}$ to $v_{t_2}$ is
\[
\Delta^{(n)}_i\;=\;E^{(n)}(t_1)+E^{(n)}(t_2)-2\inf_{[t_1,t_2]}E^{(n)}
\;=\;\sqrt{2n}\;\Delta_{m,i}\,e^{(n)},
\]
with the notation
\[\Delta_{m,i} f=\left\vert f
\parth{\frac {i-1}{m}}+f
\parth{\frac {i}{m}}-2\inf_{[\frac{i-1}{m},\frac{i}m]}f
\right \vert.
\]
The walk $(E^{(n)}(t))_{t_1\leq t\leq t_2}$ is decomposed along the
branch from $v_{t_1}$ to $v_{t_2}$ into a sequence
\[
E_0,\epsilon_1,E_1,\epsilon_2,\ldots,E_{k-1},\epsilon_{k},E_k,
\]
where the $\Delta^{(n)}_i=k$ passages on the branch separate
$\Delta^{(n)}_i$ subtrees (each coded by a $E_i$) of total size
(number of edges) $N^{(n)}_i=(2n/m-\Delta^{(n)}_i)/2$.

In this decomposition, conditionally given $\Delta^{(n)}_i=k$ (and
$N^{(n)}_i=n'=n/m-k/2$), the forest satisfies
\[
[(E_0\cdots E_{k},V_0\cdots V_{k})]\;
\mathop{=}^{\textrm{law}}\;
(E^{k+1,n'},V^{k+1,n'}),
\]
upon resetting all walks to start at zero. Under the same conditions,
the variation of labels on edges of $B$ and $C$ are i.i.d. random
variables $\zeta(\epsilon)$, uniform on $\{-1,0,+1\}$: the random
variable
\[
\bar H^{(n)}_i\;=\;\sup_{0\leq \ell\leq
\Delta^{(n)}_i}\;\sum_{j=0}^{\ell} \zeta(\epsilon_j)
\]
has tail distribution bounded by Azuma's inequality.

Our strategy is to bound $\Delta^{(n)}_i$ and,
conditionally given $\Delta^{(n)}_i$,
to bound separately the variations on the forest and on the branch:
\begin{itemize}
\item
First fix $\alpha>0$ and let $a_1$ denote the tail probability for
$\Delta^{(n)}_i$:
\begin{eqnarray*}
a_1&=&\mathbb{P}\parth{\Delta^{(n)}_i\geq m^{-\alpha} \sqrt{2n}}
\;=\;\mathbb{P}\parth{\Delta_{m,i}e^{(n)}\geq m^{-\alpha}}.
\end{eqnarray*}
\item On the complementary event for $\Delta^{(n)}_i$ we consider the
variation in a forest,
\begin{eqnarray*}
a_2&=&\mathbb{P}\parth{\sup_{0\leq t\leq 2N^{(n)}_i}
\left|V^{\Delta^{(n)}_i+1,N^{(n)}_i}(t)\right|\,\ge\,\frac{\delta'}4\, n^{1/4}\
\textrm{ and }\;\Delta_{m,i} e^{(n)}\le m^{-\alpha}},
\\
&\leq&2\mathbb{P}\parth{\sup_{0\leq
t\leq2N^{(n)}_i}{V^{\Delta^{(n)}_i+1,N^{(n)}_i}(t)}\,\ge\,\frac{\delta'}4\,
n^{1/4}\
\textrm{ and }\;\Delta_{m,i} e^{(n)}\le m^{-\alpha}},
\end{eqnarray*}
\item
and the variations on the branch from $t_1$ to $t_2$,
\begin{eqnarray*}
a_3&=&\mathbb{P}\parth{
\sup_{0\leq \ell\leq \Delta^{(n)}_i}\;\Big|\smash{\sum_{j=0}^{\ell}}
\zeta(\epsilon_j)\Big|\;\ge\;\frac{\delta'}4\,n^{1/4}\;
\textrm{ and }\;\Delta^{(n)}_i\le
m^{-\alpha}n^{1/4}}, \\
&\leq&2\mathbb{P}\parth{\bar H_i^{(n)}\,\ge\,\frac{\delta'}4\,n^{1/4}\;
\textrm{ and }\;\Delta^{(n)}_i\le m^{-\alpha}n^{1/4}}.
\end{eqnarray*}
\end{itemize}
In view of the decomposition we have
\[p_{m,i}(n)\; \le \;a_1+a_2+a_3.\]

We now bound separately each term $a_1$, $a_2$, and $a_3$.
Let us start with $a_2$ and write
\[a_2\;=\;\sum_{k=0}^{m^{-\alpha}\sqrt {2n}}\mathbb{P}\parth{
\sup_{0\leq t\leq2n'}{V^{k+1,n'}(t)}\,
\ge\,\frac{\delta'}4\,n^{1/4}\
}\mathbb{P}\parth{\Delta^{(n)}_i=k},\]
where $n'=n/m-k/2$.
In order to apply the tail estimate of the previous subsection
(Proposition~\ref{pro:loosebound}), we need
\[y\;:=\;\frac{\delta'}4(n/n')^{1/4}\;\ge\;y_0.\]
Since $n'\le n/m$ this condition is satisfied
as soon as
\[
\frac{\delta'}4\ m^{1/4}\ \ge y_0,
\quad\textrm{ that is }\quad
m\ge(4y_0/\delta')^4.\]
Then,
\begin{eqnarray*}
\mathbb{P}\parth{
\sup_{0\leq t\leq2n'}{V^{k+1,n'}(t)}
\,\ge\,\frac{\delta'}4\,n^{1/4}\ }
\;\le\; \exp\parth{-\delta' m^{1/4}/4}.
\end{eqnarray*}
The latter bound being independent of $k$ we obtain
\[
a_2\;\le\; 2\exp\parth{-\delta' m^{1/4}/4}.
\]

Let us now turn to $a_3$. Conditionally given that $\Delta^{(n)}_i=
k$, the maximum on the branch $\bar H_i^{(n)}$ is distributed as the
maximum of a random walk with $k$ steps that are independent
increments uniform on $\{-1,0,1\}$. Thus, applying Azuma's inequality
and the reflection principle, we have, for $k\le m^{-\alpha}\sqrt
{2n}$,
\[
\mathbb{P}\parth{
\left.\bar H_i^{(n)}
\,\ge\,\frac{\delta'}4\,n^{1/4}\ \right\vert
\ \Delta^{(n)}_i= k}
\;\le\; 2\, \exp\parth{-\frac{{{\delta}'}^2}{32\sqrt 2}\ m^{\alpha}}.
\]
Again the latter bound is independent of $k$ so that
\[
a_3\
\le\ \exp\parth{-\frac{{{\delta}'}^2}{32\sqrt 2}\ m^{\alpha}}.
\]

Finally let us deal with $a_1$. For $n$ large enough, the convergence
of $e^{(n)}$ to the normalised excursion $e$ entails
\[\left\vert a_1\ -\ \mathbb{P}\parth{\Delta_{m,i} e\ge
m^{-\alpha}}\right\vert\ \le\ \frac{\varepsilon }{2m}.
\]
Let us thus consider, for $1\leq i\leq m$, the probability $ \pi_{m,i}=
\mathbb{P}\parth{\Delta_{m,i}\, e\ge m^{-\alpha}}$, and restrict the
choice of $\alpha$ to $0<\alpha<1/2$.  The finite dimensional
distribution of the value of $e$ at two points and the minimum between
them was considered in Section~\ref{sec:fidis} (take $p=2$ in the
continuum limit of Proposition~\ref{pro:fidis} or see
\cite{SERLET}). For $i\neq 1,m$ this entails,
\[\pi_{m,i}=\frac{m^{3/2}}{\pi(\tau_1\tau_2')^{3/2}}\int_{
\begin{subarray}{c}
x>0,y>0,\\
\beta>m^{-\alpha},\
y+x>\beta>\vert x-y\vert
\end{subarray} }
xy\beta\ e^{-y^2/2\tau_1\ -x^2/2\tau_2'}
e^{-m\beta^2/2}dx dy d\beta,
\]
in which
\[
\tau_1= \frac{i-1}m, \qquad
\tau_2= \frac{i}m, \qquad\textrm{and}\quad \tau'_2=1-\tau_2.
\]
Thus
\begin{eqnarray*}
\pi_{m,i}&\le&\frac{m^{3/2}}{\pi(\tau_1\tau_2')^{3/2}}\int_{
\begin{subarray}{c}
x>0,y>0,\\
\beta>m^{-\alpha}
\end{subarray} }
xy\beta\ e^{-y^2/2\tau_1\ -x^2/2\tau_2'}
e^{-m\beta^2/2}dx dy d\beta
\\
&\le&\frac{m^{1/2}}{\pi(\tau_1\tau_2')^{1/2}}\ e^{-m^{1-2\alpha}/2}\\
&\le&\frac{m^{3/2}}{\pi}\ e^{-m^{1-2\alpha}/2}.
\end{eqnarray*}
For the remaining two cases $i=1$ and $i=m$, we have with
$\tau_1=\frac{m-1}m$,
\[
\pi_{m,1}=\pi_{m,m}=\sqrt{\frac{2m^3}{\pi \tau_1^3}}\int y^2
e^{-m\,y^2/2} e^{-y^2/2\tau_1}
1_{y>m^{-\alpha}} dy,\]
Now $x\rightarrow xe^{-x^2/a}$  is bounded by $\sqrt{\frac a{2e}}$,
so that
\begin{eqnarray*}
\pi_{m,1}
&\le& \sqrt{\frac{2m}{\pi \tau_1^2 e}}\int u e^{-u^2/2}
\ 1_{u>m^{0.5-\alpha}}du,\\
&\le& m^{1/2}\ e^{-m^{1-2\alpha}/2}.
\end{eqnarray*}
Thus for all $i$, $\pi_{m,i}$ is bounded by an exponentially
descreasing function of $m$.

The proof of the lemma is then concluded by summing the contribution of
     $a_1$, $a_2$ and $a_3$: for all $m$ there exists $n_0$ such that for
     all $n>n_0$,
\[\sum_{i=1}^m\ p_{m,i}(n)\le \frac{\varepsilon }{2}
     +C_\delta(m),
\]
in which $C_\delta(m)$ is exponentially small in a power of $m$.
\end{proof}

Taking $m$ and then $n$ large enough the tightness is proved,
and together with Lemma~\ref{pro:fidis} this concludes the proof
of Theorem~\ref{thm:converge}.

\subsection{Convergence of the profile}\label{sec:profile}
In view of the Skorohod representation theorem
     \cite[II.86.1]{ROW}, we may assume the joint existence,
on some probabilistic triple $(\Omega,A,\mathbb P)$,
     of a sequence of copies of  $X^{(n)}$,
     and of a copy of $X$ (and we keep the same notation
    $X^{(n)}$
     and  $X$, as for the original), such that,
for almost any $\omega \in \Omega$, $\parth{X^{(n)}_t(\omega)}_{0\le t\le 1}$
     converges to
     $\parth{X_t(\omega)}_{0\le t\le 1}$ in the Skorohod topology
     of  $D([0,1],\mathbb{R}^2)$. In this section we build copies of
$F_n$ and $F$, such
that, almost surely, $F_n$  converges to $F$.

First, the set \[\Omega_1=\left\{\omega\left\vert
\ t\longrightarrow X_t(\omega)\textrm{ is continuous}\right.\right\}\]
has probability 1, so that   uniform convergence  of
$\parth{X^{(n)}_t(\omega)}_{0\le t\le 1}$ to
$\parth{X_t(\omega)}_{0\le t\le 1}$  holds almost surely.  Set
\begin{eqnarray*}
W^{(n)}_{\min}&=&\inf_{0\le s\le 1}\
W^{(n)}_s(\omega),\\
\delta_n( \omega)&=&\sup_{0\le s\le 1}\left\vert
W^{(n)}_s(\omega)-\hat W_s(\omega)\right\vert,\\
\Psi_n(x, \omega)&=&\int_0^1 \mathbf{1}_{\hat W^{(n)}_s(\omega)\le x}\ ds,\\
\tilde F_n(x)&=&\Psi_n\parth{W^{(n)}_{\min}+x, \omega},\\
\Psi(x, \omega)&=&\int_0^1 \mathbf{1}_{\hat W_s(\omega)\le x}\ ds\\
&=&\mathcal J\parth{(-\infty,x]}.
\end{eqnarray*}
It
follows from general results on superprocessus that
$\Phi$, the distribution function of the random measure ISE,
is almost surely continuous \cite{LeGallPerso}.
Now from
\[\Psi(x-\delta_n)\le\Psi_n(x)\le\Psi(x+\delta_n)\]
and the almost sure continuity of $\Psi$, it follows that the set
\[\Omega_3=\{\omega\ \vert\ \forall  x,
\ \ \lim_n\Psi_n(x, \omega)=\Psi(x, \omega)\textrm{ and }
x\rightarrow\Psi(x, \omega)
\textrm{ is continuous}\}\]
has probability 1.  In $\Omega_3$, as
we deal with increasing functions,
\textit{uniform} convergence of
$\Psi_n$ to $\Psi$ holds true. Hence,
on the set
\[\Omega_4=\{\omega\in\Omega_3\ \vert\
\ \ \lim_nW^{(n)}_{\min}(\omega)=W_{\min}(\omega)\},\]
i.e. almost surely, uniform convergence of
$\tilde F_n$ to $F$ holds true.

On the other hand, we explain below why $\tilde F_n$
is close to some copy of $F_n$, that we shall denote
$F_n$ too. As in section \ref{sec:contourpair},
from $X^{(n)}$ one recovers
a random uniform   contour pair $(E^{(n)},V^{(n)})$,
and a random uniform
embedded tree $U_n\in\bar{\mathcal{U}}_n$.
Now we can choose at random
a well labelled tree $W_n$ in the
conjugacy class of $U_n$ in such a way  that
Theorem~\ref{thm:coupling}  holds for $(W_n,U_n)$.
Theorem~\ref{thm:coupling}
entails that $\widehat \lambda^{(n)}$ and
   $\widehat \Lambda^{(n)}$ have the same asymptotic
behavior, and $F_n$ is just  $\widehat \lambda^{(n)}$ suitably normalised.
So we have now to establish the relation between
   $\widehat \Lambda^{(n)}$ and $\tilde F_n(x)$. First, set
\[\widehat f^{(n)}_y=2n\ \tilde F_n\parth{y(8n/9)^{-1/4}, \omega}.\]
As we have
\[m_n=(8n/9)^{1/4}\ W^{(n)}_{\min},\]
$\widehat f^{(n)}_y$  is the number of visits of the contour traversal of
$U_n$ to a node whose label is not larger than $m_n+y$, but the
number of visits,
   by the contour traversal, of a given node, is exactly 1 plus the
number of children
of this node, so that
\[\widehat \Lambda^{(n)}_y+\widehat \Lambda^{(n)}_{y-1}-1\le
\widehat f^{(n)}_y\le\widehat \Lambda^{(n)}_y+\widehat \Lambda^{(n)}_{y+1}.\]
Hence, due to
Theorem~\ref{thm:coupling},
\[2\widehat \lambda^{(n)}_{y-3}-1\le
\widehat f^{(n)}_y\le 2\widehat \lambda^{(n)}_{y+3},\]
or, equivalently,
\[\frac n{n+1}\ \tilde F_n\parth{x-cn^{-1/4}, \omega}\le
F_n\parth{x, \omega}\le\frac n{n+1}\ \tilde F_n\parth{x+cn^{-1/4}, \omega}
+\frac 1{2n+2},\]
$c$ being a constant. That is, on $\Omega_4$,
   uniform convergence of $F_n$ to $F$ holds true.

\subsection*{Acknowledgements} We thank David Aldous,
Jean-Fran\c{c}ois Le Gall and
Balint Virag for stimulating discussions.

\bibliography{resume}

\end{document}